\documentclass[leqno]{amsart}
 \usepackage{amsaddr}
\usepackage[a4paper,margin=1.5in]{geometry}
\usepackage{mathtools,amsmath,amsthm,amsfonts,amssymb,enumerate}
\usepackage[T1]{fontenc}
\usepackage[final, pdfpagelabels, pdfstartview = {FitH}, bookmarks, colorlinks=true, plainpages = false, linktoc=all, linkcolor=red, citecolor=blue, urlcolor=blue,filecolor=black]{hyperref}
\mathtoolsset{showonlyrefs=true}
\usepackage{url}
\usepackage{palatino}
\usepackage{euscript}
\usepackage{xcolor}
\usepackage{enumitem}
\usepackage{pdfpages} 
\usepackage{pdfrender, tikz-cd, rotating}
\usepackage{subcaption}
\makeatletter
\def\mathcolor#1#{\@mathcolor{#1}}
\def\@mathcolor#1#2#3{%
  \protect\leavevmode
  \begingroup
    \color#1{#2}#3%
  \endgroup
}
\makeatother

\usepackage{pstricks}
\usepackage[all]{xy}
\SelectTips{cm}{}
\usepackage{pinlabel}

\newtheorem{thm}{Theorem}[subsection]
\newtheorem{lem}[thm]{Lemma}

\newtheorem{prop}[thm]{Proposition}

\theoremstyle{definition}
\newtheorem{defn}[thm]{Definition}

\theoremstyle{remark}
\newtheorem{rem}[thm]{Remark}

\numberwithin{thm}{section}

\newcommand{\nc}{\newcommand} 
\nc{\mb}{\mathbb}
\nc{\mc}{\mathcal}
\nc{\mf}{\mathfrak}

\nc{\ic}{\mathbf{IC}}

\newcommand{\bbN}{\mathbb{N}}

\newcommand{\leqH}{\geq_{\mathbf{H}}}

\newcommand{\bfH}{\mathbf{H}}

\newcommand{\undH}{\underline{\bfH}}

\usepackage{}

\DeclareMathOperator{\aut}{Aut}

\newcommand{\joins}[5]{\mathcal{J}_{#1,#2}^{[#3,#4]}(#5) }

\definecolor{Pink}{rgb}{1.,0.75,0.8}
\definecolor{mygray1}{gray}{0.20}
\definecolor{mygray2}{gray}{0.40}
\definecolor{mygray3}{gray}{0.60}
\definecolor{mygray4}{gray}{0.80}

\definecolor{ao(english)}{rgb}{0.0, 0.5, 0.0}

\usepackage{mathrsfs}

\usepackage{bbm}

\def\1{\mathbbm{1}}

\def\Today{\ifcase\month\or January\or February\or March\or
  April\or May\or June\or July\or August\or September\or
  October\or November\or December\fi\space\number\year}

\definecolor{burgundy}{rgb}{0.5, 0.0, 0.13}

\renewcommand{\hat}{\widehat}

\usepackage{bbm}

\def\1{\mathbbm{1}}

\numberwithin{equation}{section}

\newcommand{\NN}{\mathbf{N}}
\newcommand{\HH}{\mathbf{H}}
\newcommand{\MM}{\mathbf{M}}

\tikzcdset{ cells={font=\everymath\expandafter{\the\everymath\displaystyle}}, }

\author{Gaston Burrull}
\address{Faculty of Science, The University of Sydney, Carslaw Building, Eastern Ave., Camperdown Sydney, New South Wales 2006, Australia}
\thanks{The first author was supported by Beca Chile Doctorado 2017 Folio No. 72180433.}

\author{Nicolas Libedinsky}
\address{Universidad de Chile, Mathematics Las Palmeras, 3425 Casilla, 653 Santiago, Ñuñoa, 7800003, Chile}
\thanks{The second author was partially supported by FONDECYT-ANID grant 1200061.}

\author{David Plaza}
\address{Universidad de Talca, Instituto de Matemáticas, Avenida Lircay s/n Talca, 3460000, Chile}
\thanks{The third author was partially supported by FONDECYT-ANID grant 1200341.}

\title{Combinatorial invariance conjecture for $\widetilde{A}_2$}

\begin{document}

\maketitle

\begin{abstract} 
The combinatorial invariance conjecture (due independently to G. Lusztig and M. Dyer) predicts that if $[x,y]$ and $[x',y']$ are isomorphic Bruhat posets (of possibly different Coxeter systems), then the corresponding Kazhdan-Lusztig polynomials are equal, that is, $P_{x,y}(q)=P_{x',y'}(q)$. We prove this conjecture for the affine Weyl group of type $\widetilde{A}_2$. This is the first infinite group with non-trivial Kazhdan-Lusztig polynomials where the conjecture is proved.
\end{abstract}
\vskip 1cm

\section{Introduction}

Kazhdan-Lusztig polynomials $P_{x,y}(q)$ (in particular for affine Weyl groups) are of fundamental importance in the representation theory of Lie theoretic objects and in the topology and geometry of Schubert varieties. %Furthermore, the fact that their coefficients are non-negative (see \cite{EWhodge}) makes them interesting combinatorial objects. 
There is a beautiful combinatorial conjecture involving them that was predicted independently by G. Lusztig (unpublished, see \cite{Bre02}) and M. Dyer \cite{dyer1987hecke}. If this conjecture was verified it would be very surprising both from the algebraic and from the geometric perspective. 

On the algebraic side, the recursive algorithm producing Kazhdan-Lusztig polynomials in the Hecke algebra does seem to care about the specifics of the elements involved and not just about their Bruhat order. We also lack of any heuristic on why this conjecture should be true. 

On the geometric side, if $G$ is a complex Kac-Moody group, $B$ a Borel subgroup and $W$ is the associated  Weyl group, one has the Bruhat stratification of the generalized flag variety
\begin{equation*}
G/B=\biguplus_{w\in W}BwB/B.
\end{equation*}
The locally closed subvariety $X_w:=BwB/B$ is called the \emph{Schubert cell} and its closure $\overline{X}_w$ is called the \emph{Schubert variety} of $w\in W$. In this situation, the coefficients of the Kazhdan-Lusztig polynomial $P_{w,v}(q)$ are given by the dimensions of the local intersection cohomology (with middle perversity and associated with the Bruhat stratification) of $\overline{X}_v$ at any point of $X_w$ (cf. \cite[Theorem 12.2.9]{kum02}). The key point here is that $w\leq v$ in the Bruhat order is equivalent to  $\overline{X}_w\subset \overline{X}_v.$ So, if the conjecture was true, the poset defining the stratification would give the dimensions of the corresponding local intersection cohomology spaces. In similar situations this is not true, for example in the case of the nilpotent cone or in the case of isolated singularities (e.g. the cone over a smooth projective variety \cite[Proposition 2.4.4]{Ara02}), etc. 

On the other hand, there are several results that tilt the situation towards believing that the conjecture might be true. The most important in our opinion is the main result in \cite{BCM06} (following a similar but weaker result in \cite{Du03}) by F. Brenti, F. Caselli, and M. Marietti, where the conjecture is proved for any pair of Coxeter systems for lower intervals. In the notation of the abstract, lower intervals are those for which $x$ and $x'$ are the identity elements of the corresponding groups. A more recent result by L. Patimo \cite{Patimo2021} shows that the coefficient of $q$ of the Kazhdan-Lusztig polynomial in finite type ADE is invariant under poset isomorphisms. Some other interesting papers related to this conjecture are \cite{Dy91}, \cite{Br94}, \cite{Br97}, \cite{Bre02}, \cite{Br04}, \cite{incitticombinatorial}, \cite{incitti2007more}, \cite{Ma06}, and \cite{Ma18}.

In this paper we prove the conjecture when both Coxeter groups involved are the affine Weyl group of type $\widetilde{A}_2$. As said in the abstract, this is the first infinite group with non-trivial Kazhdan-Lusztig polynomials where the conjecture is proved. This is interesting new evidence in favor of the conjecture because this is the first time that the conjecture is proved for arbitrarily large posets that are not lower and have non-trivial Kazhdan-Lusztig polynomials.

We heavily use the explicit formulas for Kazhdan-Lusztig polynomials in $\widetilde{A}_2$ recently found in \cite{libedinsky2020affine} and the fact that the conjecture is known for intervals $[x,y]$ when $\ell(y)-\ell(x)\leq 4$ \cite[Exercises 5.7 and 5.8]{BB05}. The main technical ingredient in the proof are the sets
\begin{equation*}
 \{ z\in [x,y] \, |\, \ell(y) -\ell(z)=m \mbox{ and } P_{z,y}(q)=1+q  \}.
\end{equation*}  
For $1\leq m\leq 4 $ we check that they are preserved under poset isomorphisms (for a precise statement see Lemma \ref{lemma Z to Z}). It is interesting to notice that this invariant includes Kazhdan-Lusztig polynomials in its definition. Combining this and the explicit formulas mentioned above, the result follows by induction. 

It has  been known from the dawn of the theory that Bruhat intervals and  Kazhdan-Lusztig combinatorics for $\widetilde{A}_1$ are trivial: modulo isomorphism, there is one Bruhat interval of length $n\in \mathbb{N}$ and Kazhdan-Lusztig polynomials are all $1$. The  $\widetilde{A}_2$ case is radically different. We were not able to solve the classification problem of Bruhat intervals in this case (this was our first approach towards the result in this paper) due to its high complexity. As an illustration of this complexity, the number of non-isomorphic Bruhat intervals of length $n\in \mathbb{N}$ is an unbounded function on $n.$ On the other hand, Kazhdan-Lusztig polynomials in the $\widetilde{A}_2$ case were found only in 2020 \cite{libedinsky2020affine}. Similar formulas \cite{LPP21} have been found in 2021 by Patimo, the second and the third authors, in the cases of $\widetilde{A}_3$,  $\widetilde{A}_4$ (and $\widetilde{A}_5$, unpublished) for the maximal elements in their double cosets (which is called in this paper the set $\Theta$). In another unpublished work by the same authors, using geometric Satake they find similar formulas, also up to $\widetilde{A}_5$ for the lowest double Kazhdan-Lusztig cell (regions $\Theta, \Theta_1$, and $\Theta_2$ in this paper). In conclusion, the combinatorics of $\widetilde{A}_2$ has a very similar taste to that of higher ranks.

So in principle, to replicate the results of this paper in higher ranks, the key result that one would need (and that we do not have at present) is a proof of the conjecture  for intervals when $\ell(y)-\ell(x)$ is small. For example, in the case of  $\widetilde{A}_3$, we would need to know that the conjecture is valid when $\ell(y)-\ell(x)\leq 6$.

%{\red{Added by David. Esto es lo que queremos decir en formulas} } 
%Let $a_m^n$ be the number (up to isomorphism) of intervals of length $m$ in an affine Weyl group of type $\tilde{A}_{n}$. We have
%\begin{equation}
%    \lim_{m\rightarrow \infty} a_{m}^1 =1  \qquad \lim_{m\rightarrow \infty} a_{m}^2 =\infty.  
%\end{equation}
%Por ejemplo, es facil ver que $m \leq a_{2m}^2$ usando las formulas del tamano de los intervalos theta-theta que alguna vez teniamos. De hecho si $y=\theta(2m,2m)$ hay al menos $m$ intervalos no isomorfos de la forma $[x,y]$ con $\ell(y)- \ell(x)=2m$. 
\subsection{Structure of the paper} In Section \ref{preliminaries} we recall some background material about the  affine Weyl group $\tilde{A}_2$ and  some results from \cite{libedinsky2020affine}. In Section \ref{section KL formulas} we write down and prove even more explicit formulas for Kazhdan-Lusztig polynomials than the ones given in \cite{libedinsky2020affine}. In Section \ref{inv} we define two poset invariants that are useful to discard possible isomorphisms of posets and find some properties satisfied by them. Finally, in Section \ref{proof} we prove the main theorem. 

\section{Preliminaries}\label{preliminaries}
%In this section we will review some basic properties about the $\tilde{A}_2$ affine Weyl group and its Hecke algebra.

\subsection{The affine Weyl group of type \texorpdfstring{$\tilde{A}_2$}{A2 tilde}} \label{prelim A2 Weyl}
Let $W$ be the affine Weyl group of type $\tilde{A}_2$.  It is a Coxeter system generated by the simple reflections  $S=\{s_0,s_1,s_2\}$ with relations $s_i^2=\operatorname{id}$ for $i\in \{0,1,2\}$ and $(s_is_j)^3 =\operatorname{id}$ for $i\neq j \in \{0,1,2\}$. If no confusion is possible, we will sometimes denote the generators $s_0$, $s_1$, and $s_2$ by $0$, $1$, and $2$, respectively. We will also use ``label$\mod 3$" notation. For example, $135678$ stands for $s_1s_0s_2s_0s_1s_2$. As usual, we denote by $\ell (\cdot )$ and $\leq$ the length and the Bruhat order on $W$, respectively. 

The Dynkin diagram of $W$ has six symmetries (it is an ``equilateral triangle''). Each one of them induces an automorphism of $W$. We denote by $\rho\colon S\rightarrow S$ the map given  by  $\rho (s_0) = s_1$, $\rho (s_1)=s_2$ and $\rho (s_2)=s_0$.  Similarly, we consider $\sigma\colon S\rightarrow S$ the map that fixes $s_0$ and permutes $s_1$ and $s_2$.  The maps  $\rho $ and $\sigma$ extend to automorphisms of $W$, which we denote by the same symbols. We denote by $D_3$ the subgroup of $\aut (W)$ generated by $\rho$ and $\sigma$. We denote by $\iota\colon W \rightarrow W$ the inversion anti-automorphism, that sends
$x\mapsto x^{-1}$.  We define $G$ to be the subgroup of $\mbox{Sym}(W)$ generated by $\rho$, $\sigma$, and $\iota$. 

Let $\underline{y} = (r_1, r_2, \ldots,r_k)$ be an expression for $y$ (i.e., $r_i\in S$ for all $i$). If $1<i<k$, we say that there is a \emph{braid triplet in position $i$} if $(r_{i-1}, r_i, r_{i+1}) = (s_j,s_{j'},s_j)$ with $j\neq j'$. We define 
the \emph{distance} between a braid triplet in position $i$ and a braid triplet in position $j>i$ to be the number $ i-j-1$.
The following is \cite[Lemma 1.1]{libedinsky2020affine}. 

\begin{lem}\label{lemma braid triplets}
An expression $\underline{w}$ without adjacent simple reflections is reduced if and only if the distance between any two braid triplets is odd.
\end{lem} 

For any positive integer $n$, we define $x_n\coloneqq 123\cdots n \in W$. Note this is a reduced expression  by Lemma \ref{lemma braid triplets} (there are no braid triplets), therefore $\ell(x_n)=n$. Let us define
\begin{equation}
    X \coloneqq \bigcup_{\tau \in G} \tau (\{x_n \, |\, n\geq 1\}). 
\end{equation}
For any pair $(m,n)$ of non-negative integers, we define 
\begin{equation} \label{def thetas}
\theta(m,n)\coloneqq 1234\cdots (2m+1)(2m+2)(2m+1)\cdots(2m-2n+1).
\end{equation}
This is a reduced expression of $\theta(m,n)$ by Lemma \ref{lemma braid triplets} (there is only one braid triplet). In particular, $\ell(\theta(m,n))=2m+2n+3$. It is easy to see that $\theta(m-1,n-1)< \theta(m,n)$. We define
\begin{equation}
    \Theta \coloneqq \bigcup_{\tau \in G} \tau (\{\theta(m,n) \, |\, m,n\geq 0 \}). 
\end{equation}
The elements in $\Theta$ can be characterized as those with left and right descent set formed  exactly by two simple reflections. Thus, for any $\theta(m,n)$ there exists a unique $s_{m,n}\in S$ such that $$\ell (\theta(m,n)) < \ell (\theta(m,n)s_{m,n} ).$$ On the other hand, $s_0$ is the unique simple reflection that is not in the left descent set of any $\theta(m,n)$. We define
\begin{equation}
 \Theta_1  \coloneqq \bigcup_{\tau \in G  } \tau (\{\theta(m,n)s_{m,n} \, |\, m,n\geq 0 \}) \quad \mbox{and} \quad    \Theta_2  \coloneqq \bigcup_{\tau \in G } \tau (\{s_0\theta(m,n)s_{m,n} \, |\, m,n\geq 0 \}).
\end{equation}
 Note that $\ell(\theta(m,n)s_{m,n})=2m+2n+4$ and $\ell(s_0\theta(m,n)s_{m,n})=2m+2n+5$.\\

% \noindent Let $y=s_0\theta(m,n)s, $ where $s=s_{m,n}$. It follows from Equation \eqref{def thetas} and Lemma \ref{lemma braid triplets} that $y$ admits the following reduced expression
% \begin{equation}\label{rex expression y in remark}
%   \underline{y} \coloneqq  01234\cdots (2m+1)(2m+2)(2m+1)\cdots (2m-2n+1)(2m-2n). 
% \end{equation}
% Using \cite[Lemma 1.1]{libedinsky2020affine} we can see that the unique simple reflections that can be removed from $\underline{y}$ in order to obtain a reduced expression are the first two, the last two or the two simple reflections denoted by $(2m+1)$ in \eqref{rex expression y in remark}.

Given $y\in W$ we define $C_y:=\{ w\in W \, |\, w <y \mbox{ and } \ell(w)=\ell(y)-1 \}$. 

%The following lemma follows from Lemma \ref{lemma braid triplets}.

\begin{lem}\label{corolario ganador}
 Let $m$ and $n$ be positive integers. Set $s=s_{m,n}$.  If $y=s_0\theta(m,n)s$ then 
%  \begin{equation}
%     C_y=\{ z_1\rho^2(s)s, z_2\rho(s)s, ys, s_0s_1s_0y, s_0s_2s_0 y, s_0y\}.
% \end{equation}  
 \begin{gather}
    C_y=\{\theta(m,n)s, \rho (\theta(m-1,n+1))s, s_0\theta(m-1,n+1),  \rho^2 (\theta(m+1,n-1))s, \\ s_0\theta(m+1,n-1), s_0\theta(m,n)   \}.
\end{gather} 
\end{lem}
 
\begin{proof}
Every element in $C_y$ can be obtained by removing a simple reflection from a reduced expression for $y$ as long as the resulting expression is reduced. By Lemma \ref{lemma braid triplets} 
 \begin{equation}\label{rex expression y in remark}
  \underline{y} \coloneqq  01234\cdots (2m+1)(2m+2)(2m+1)\cdots (2m-2n+1)(2m-2n) 
\end{equation}
is a reduced expression for $y$. Using  Lemma \ref{lemma braid triplets} once again we have that the elements of $C_y$ are obtained from $\underline{y}$ by removing any of the first two, the last two or the two simple reflections denoted by $(2m+1)$ in \eqref{rex expression y in remark}.  In formulas,
\begin{align*}
\widehat{0}123\cdots (2m+1)(2m+2)(2m+1)\cdots (2m-2n+1)(2m-2n) &= \theta(m,n)s;\\
0\widehat{1}23\cdots (2m+1)(2m+2)(2m+1)\cdots (2m-2n+1)(2m-2n) &= \rho (\theta(m-1,n+1))s;\\
0123\cdots \widehat{(2m+1)}(2m+2)(2m+1)\cdots (2m-2n+1)(2m-2n) &= s_0\theta(m-1,n+1);\\
0123\cdots (2m+1)(2m+2)\widehat{(2m+1)}\cdots (2m-2n+1)(2m-2n) &= \rho^2 (\theta(m+1,n-1))s;\\
0123\cdots (2m+1)(2m+2)(2m+1)\cdots \widehat{(2m-2n+1)}(2m-2n) &= s_0\theta(m+1,n-1);\\
0123\cdots (2m+1)(2m+2)(2m+1)\cdots (2m-2n+1)\widehat{(2m-2n)} &= s_0\theta(m,n) , 
\end{align*}
% \begin{align*}
% \widehat{0}1234\cdots (2m+1)(2m+2)(2m+1)\cdots (2m-2n+1)(2m-2n) &= s_0 y;\\
% 0\widehat{1}234\cdots (2m+1)(2m+2)(2m+1)\cdots (2m-2n+1)(2m-2n) &= s_0s_1s_0y;\\
% 01234\cdots \widehat{(2m+1)}(2m+2)(2m+1)\cdots (2m-2n+1)(2m-2n) &= z_2\rho(s)s;\\
% 01234\cdots (2m+1)(2m+2)\widehat{(2m+1)}\cdots (2m-2n+1)(2m-2n) &= s_0s_2s_0 y;\\
% 01234\cdots (2m+1)(2m+2)(2m+1)\cdots \widehat{(2m-2n+1)}(2m-2n) &= z_1\rho^2(s)s;\\
% 01234\cdots (2m+1)(2m+2)(2m+1)\cdots (2m-2n+1)\widehat{(2m-2n)} &= ys, 
% \end{align*}
where $\hat{k}$ means letter $k$ is omitted. The result follows. 
\end{proof}

\begin{figure}[ht] 
  \begin{subfigure}[t]{0.475\linewidth}
    \centering
\includegraphics[scale=0.3]{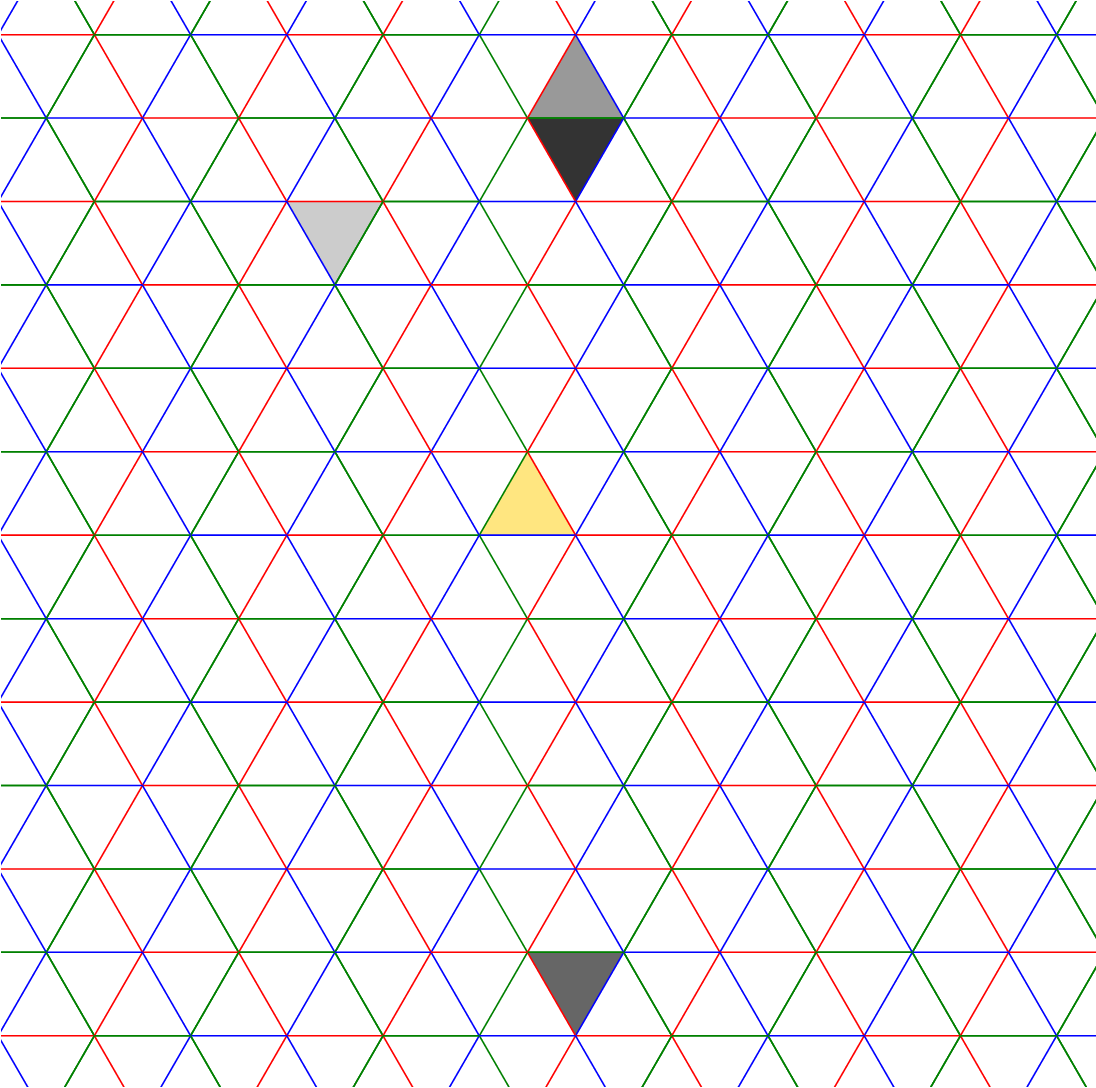}\\
    \caption{The yellow triangle corresponds to the identity. From lightest to darkest gray the other colored elements are $x_7$, $\theta(1,2)s_1$, $s_0\theta(1,2)s_1$, and $\theta(1,2)$.}
    \label{fig: W}
  \end{subfigure}\hfill%  
  \begin{subfigure}[t]{0.475\linewidth}
        \centering
\includegraphics[scale=0.43,trim={1cm 1cm 1cm 1cm},clip]{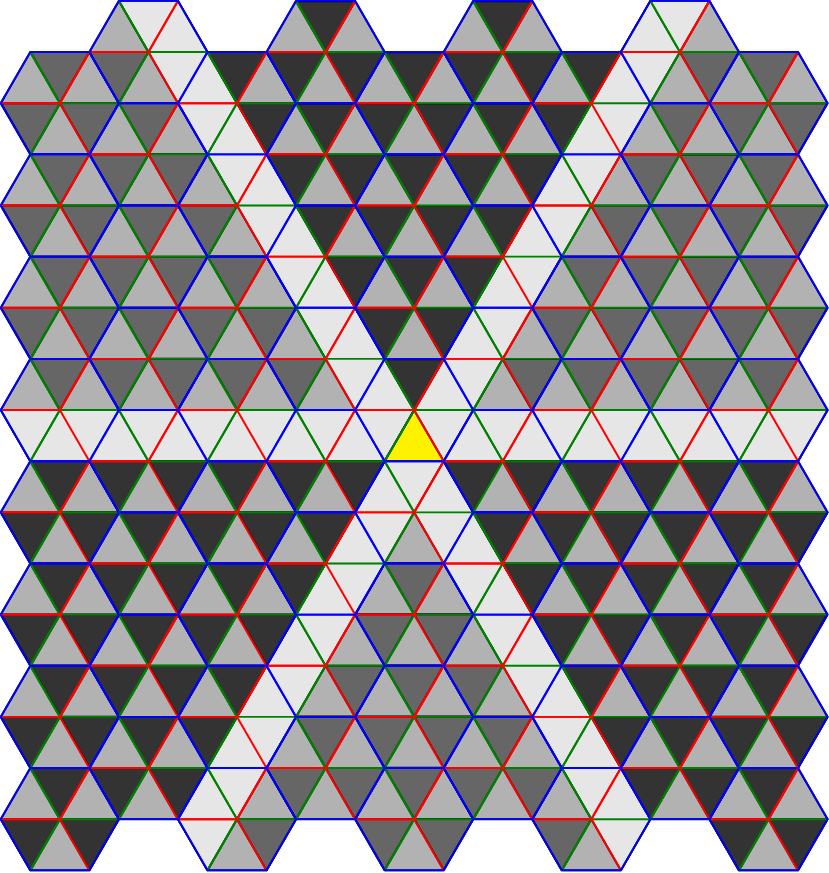}\\
    \caption{From lightest to darkest gray: $X$, $\Theta_1$, $\Theta_2$, and $\Theta$.}
    \label{fig: four regions}
  \end{subfigure}\hfill%  
  \caption{Geometric realization of $W$ and the four regions.}
\end{figure}

It will be convenient for our purposes to recall the realization of $W$ as the group of isometric (or affine) transformations of the plane generated by the three reflections with respect to the lines supporting the edges of  any  equilateral triangle of the tessellation $T$ illustrated in Figure \ref{fig: W}. Let $s_0$ (resp. $s_1$, $s_2$) be the orthogonal reflection through the  line supporting the blue (resp. green, red) edge of  the yellow  equilateral triangle. The group $W$ acts simply transitively in $T$, thus there is a bijection between elements in $ W$ and equilateral triangles in $T.$ This bijection sends $w\in W$ to the triangle $w\cdot \Delta_0,$ where $\Delta_0$ is the yellow triangle (that corresponds to the identity in $W$). If $e$ is an edge of $\Delta_0$ colored with color $c$ then in the triangle $w\cdot \Delta_0$, the edge $w\cdot e$ is also colored $c.$ Henceforth, we do not distinguish between elements of $W$ and their corresponding equilateral triangles. We also identify $s_0$, $s_1$, and $s_2$ with the colors blue, green, and red, respectively. This last identification is useful because if $\Delta$ and $\Delta'$ share an edge colored, say, green, one knows that if $w$ is an expression for $\Delta$ then $ws_1$ is an expression for $\Delta'.$

%Let us explain more algorithmically  how this bijection works. Let us identify $s_0$, $s_1$, and $s_2$ with the colors blue, green, and red, respectively. As we said, the yellow triangle corresponds to the identity. Given an element $w\in W$ and an expression $s_{i_1}s_{i_2}\ldots s_{i_k}$ (not necessarily reduced) of $w$ we will define a sequence of triangles $(\Delta_0,\Delta_1,\ldots , \Delta_k)$, each of which has its edges painted with blue, green and red in some order.  $\Delta_0$ is the identity triangle with the colors painted as in Figure \ref{fig: W}. Then $\Delta_{j+1}$ is obtained from $\Delta_j$ by reflecting it through the side colored with $s_{i_j}$, and we obtain a new coloring of the edges of $\Delta_{j+1}$ accordingly. We identify $w$ with $\Delta_{k}$. It is well-known that this  identification does not depend on the choice of an expression for $w$ (this is also true for the colors of the edges of the triangles) and that it defines a bijection between $W$ and the set of triangles in the tessellation. Henceforth, we do not distinguish between elements of $W$ and triangles. 

% \begin{figure}[ht]
%     \centering
% \includegraphics[scale=0.4]{CICA2figs/LDfig1.eps}\\
%     \caption{The yellow triangle corresponds to the identity. From lightest to darkest gray the other colored elements are $x_7$, $\theta(1,2)s_1$, $s_0\theta(1,2)s_1$, and $\theta(1,2)$.}
%     \label{fig: W}
% \end{figure}

Using the geometric realization of $W$ it is not hard to see that $ W\setminus \{\operatorname{id}\} = X \uplus \Theta \uplus \Theta_1 \uplus \Theta_2$, where $\uplus$ denotes a disjoint union, as is illustrated in Figure  \ref{fig: four regions}.

The geometric realization of $W$ allows us to describe certain lower intervals in an elegant way as the convex hull of a given set. We denote by  $\operatorname{cen}(w)$ the centroid of $w\in W$. Let $W_f\coloneqq \langle s_1,s_2\rangle \leq W$ (it is a dihedral group of order $6$).  For $w$ in $W$, let $\mathcal{C}_w$ denote the convex hull of the set $\{\operatorname{cen}(y w) \, |\, y\in W_f\}$.

The following lemma appears in the proof of \cite[Lemma 1.4]{libedinsky2020affine}. 

\begin{lem} \label{menores que theta}
Let $m$ and $n$ be non-negative integers. 
Then, 
\begin{equation*}
[\operatorname{id}, \theta(m,n)]=\{w\in W  \, |\, \operatorname{cen}(w)\in \mathcal{C}_{\theta(m,n)}\}.
\end{equation*}
\end{lem}

\begin{figure}[ht] 
  \begin{subfigure}[b]{0.475\linewidth}
    \centering
    \includegraphics[scale=0.3]{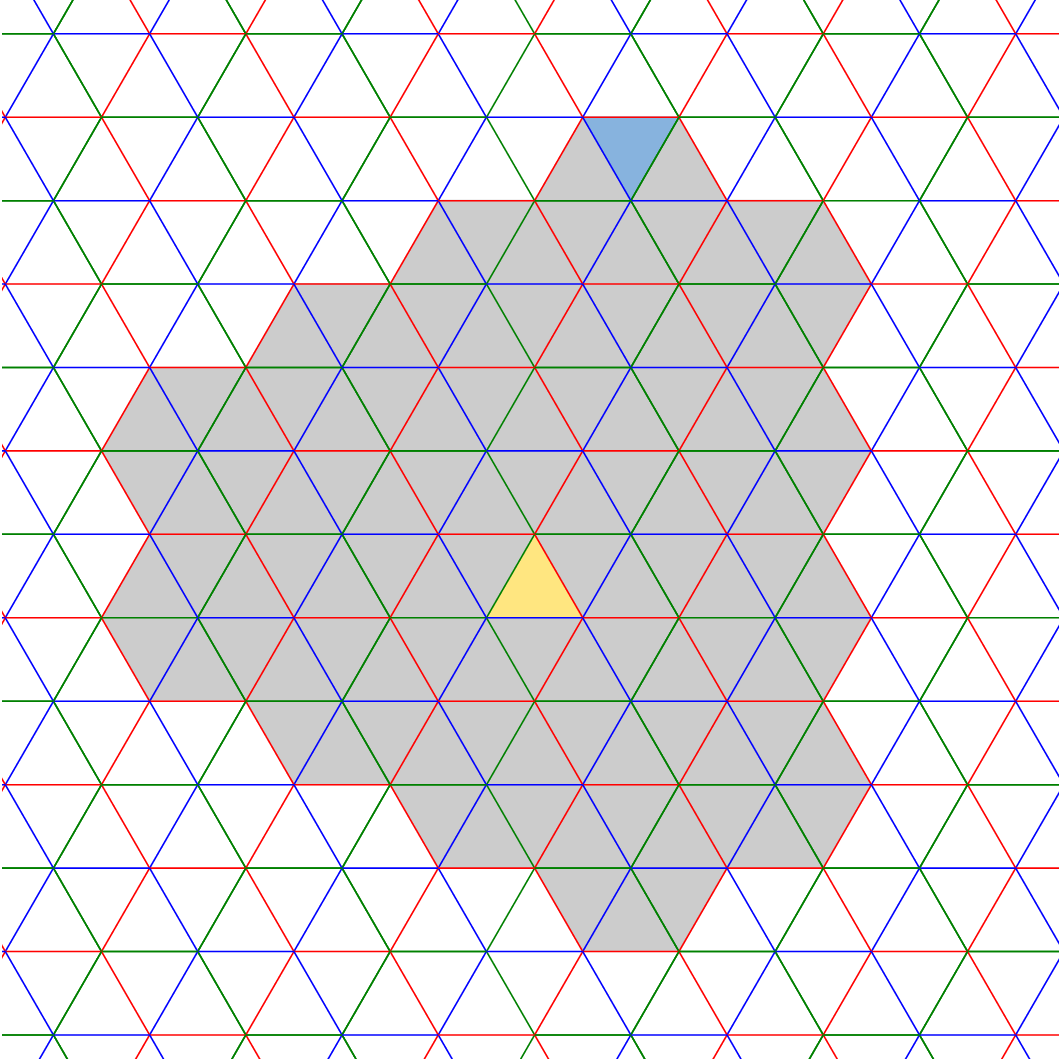}
    \caption{ Lower interval $[\operatorname{id}, \theta(1,3)]$. The blue triangle corresponds to $\theta(1,3)$. The gray triangles correspond to the non-endpoints of $[\operatorname{id}, \theta(1,3)]$.}
    \label{fig: new intervals izquierda} 
  \end{subfigure}\hfill%  
  \begin{subfigure}[b]{0.475\linewidth}
    \centering
    \includegraphics[scale=0.3]{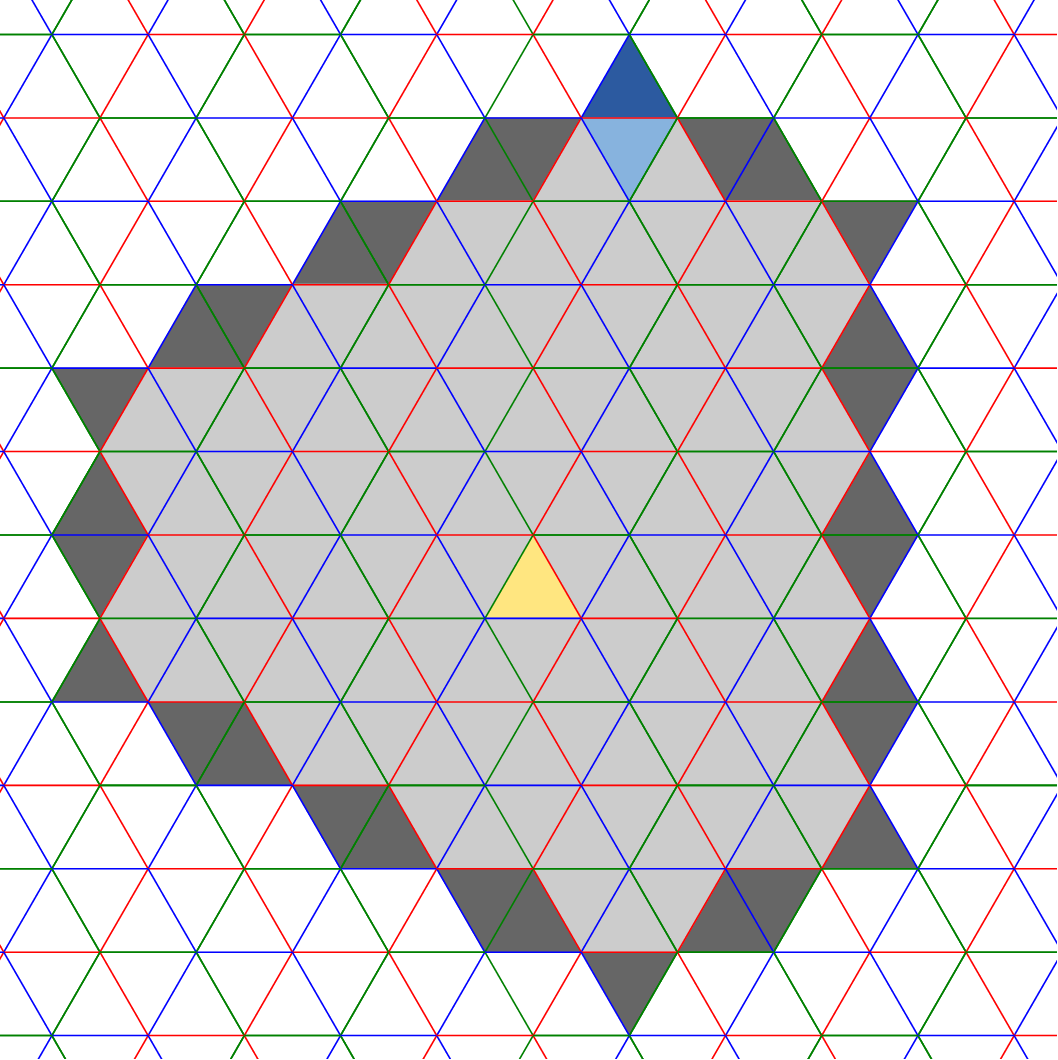}
    \caption{Lower interval $[\operatorname{id}, \theta(1,3)s_2]$. The dark blue triangle corresponds to $\theta(1,3)s_2$. The dark gray triangles together with the dark blue triangle correspond to $\partial_{1,3}s_{2}$.} 
    \label{fig: new intervals derecha} 
  \end{subfigure}\hfill%  
  \caption{Geometric description of  $[\operatorname{id}, \theta(1,3)]$ and $[\operatorname{id}, \theta(1,3)s_{2}]$.}
  \label{fig: all theta's friend} 
\end{figure}

For instance, the interval $[\operatorname{id}, \theta(1,3)]$ is illustrated in Figure \ref{fig: new intervals izquierda}. We recall that $s_1$ and $s_2$ have been identified with the reflections through the lines that support the green and red sides of the identity triangle, respectively.

% The following lemma is the generalization of Lemma \ref{menores que theta} to some of the other regions.

\begin{lem} \label{menores que s0theta}
Let $m$ and $n$ be non-negative integers. Set $s=s_{m,n}$. Then,
\begin{equation}\label{lala}
[\operatorname{id},\theta(m,n)s] = [\operatorname{id}, \theta(m,n)] \uplus \partial_{m,n}s,
\end{equation}
where $\partial_{m,n}$ is the set formed by all the elements $x \in [\operatorname{id}, \theta(m,n)]$ such that $x$ has a side belonging to the boundary of $[\operatorname{id}, \theta(m,n)]$. %In particular,
%\begin{equation*}
%[\operatorname{id}, \theta(m,n)s]=\{w\in W: \operatorname{cen}(w)\in \mathcal{C}_{\theta(m,n)s}\}.
%\end{equation*}
%Furthermore, 
%\begin{align*}
%[\operatorname{id}, s_0\theta(m,n)]&=\{w\in W: \operatorname{cen}(w)\in \mathcal{C}_{s_0\theta(m,n)}\cup s_0\mathcal{C}_{s_0\theta(m,n)}\},\\
%[\operatorname{id}, s_0\theta(m,n)s]&=\{w\in W: \operatorname{cen}(w)\in \mathcal{C}_{s_0\theta(m,n)s} \cup s_0\mathcal{C}_{s_0\theta(m,n)s}\}.
%\end{align*}
\end{lem}
\begin{proof}
It is clear that $[\operatorname{id},\theta(m,n)s] \supset [\operatorname{id}, \theta(m,n)] \uplus \partial_{m,n}s$. For the other inclusion, we first notice that if $w\leq \theta(m,n)s$ then we either have  $w\leq \theta(m,n)$ or $w= xs$ for $x\leq \theta(m,n)$. If $w\leq \theta(m,n)  $ then there is nothing to prove. Thus we can assume $w\not \leq \theta(m,n)$. Therefore  $w= xs$ for some $x\leq \theta(m,n)$. Since $x$ is inside  and $w=xs$ is outside the interval $[\operatorname{id}, \theta(m,n)]$, their shared ($s$-colored) side is in the boundary of $[\operatorname{id}, \theta(m,n)]$. It follows that  $x\in \partial_{m,n}$ and $w\in \partial_{m,n}s$.

In the proof of \cite[Lemma 1.4]{libedinsky2020affine}, the authors showed that the boundary of the interval $[\operatorname{id},\theta(m,n)]$ is colored by $s$. This implies that for $x\in \partial_{m,n}$ the element $xs$ is outside the interval $[\operatorname{id},\theta(m,n)]$, so the union in Equation \eqref{lala} is indeed disjoint. 
% It is clear that $[\operatorname{id},\theta(m,n)s] \supset [\operatorname{id}, \theta(m,n)] \uplus \partial_{m,n}s$. For the other inclusion, we first notice that if $w\leq \theta(m,n)s$ then we either have  $w\leq \theta(m,n)$ or $w= xs$ for $x\leq \theta(m,n)$. If $w\leq \theta(m,n)  $ then there is nothing to prove. We can assume $w= xs$ for $x\leq \theta(m,n)$ and $w\not \leq \theta(m,n)$. Since $x$ is inside  and $w=xs$ is outside the interval $[\operatorname{id}, \theta(m,n)]$, their shared ($s$-colored) side is in the boundary of $[\operatorname{id}, \theta(m,n)]$, therefore $x\in \partial_{m,n}$ and $w\in \partial_{m,n}s$.  
%The last two equations are proved in a similar way as the first one. The unique difference is noting the following two facts: if $w\leq s_0\theta(m,n)$ we have either $w\leq \theta(m,n)$ or $w= s_0d$ for $d\leq \theta(m,n)$.  If $w\leq s_0\theta(m,n)s$ we have either $w\leq s_0\theta(m,n)$ or $w= ds$ for $d\leq s_0\theta(m,n)$.
\end{proof}

Lemma \ref{menores que s0theta} is illustrated in Figure \ref{fig: new intervals derecha} for the pair $(m,n)=(1,3)$.

% The intervals $[\operatorname{id},\theta(1,3)s_2]$ and $[\operatorname{id},s_0\theta(1,3)s_2]$ are pictured in Figure \ref{fig: new intervals derecha} and Figure \ref{fig: lema28}, respectively.

\begin{rem}\label{algebraic description theta}
For the more algebraic or combinatorially minded readers, we give an explicit description of the lower intervals without a proof. Although one can use these descriptions in some of the results that follow, we strongly prefer the geometric versions used in the paper. 

For any pair of non-negative integers we define $a_{m,n}\coloneqq 2m+n$ and $b_{m,n}\coloneqq m+2n$.  Let  $f$ be the bijection of $\{0,1,2\}$ that fixes $0$ and permutes $1$ and $2$. The set 
\begin{equation*}
Y\coloneqq [\mbox{id}, \theta(m,n)]\setminus \{\mbox{id}\} 
\end{equation*} 
is partitioned as follows (for all the sets below, we assume that $i\in \{0,1,2\}$):
\begin{align} 
    X\cap Y   = &\,  \{  \rho^i  (x_k)\, |\, 1\leq k< a_{m,n}+3-i \} \uplus \{ \rho^i  (\sigma(x_k)) \, |\,  2\leq k< b_{m,n}+3-f(i)\} \\
        \Theta \cap Y  = &\,  \{\rho^{i}(\theta(p,q)) \, |\, a_{p,q}< a_{m,n}+1-i \mbox{ and }  b_{p,q} < b_{m,n}+1-f(i)\}\\
   \label{super description}      \Theta_1 \cap Y  = &\,   \{ \rho^{i}(\theta(p,q)s_{p,q})\, |\,  a_{p,q}< a_{m,n}-i \mbox{ and }  b_{p,q} < b_{m,n}-f(i)  \} \,  \uplus  \\
 \nonumber & \,  \{ \rho^{i}(s_0\theta(p,q)) \, |\, a_{p,q}< a_{m,n}-2 +f(i) \mbox{ and }  b_{p,q} <   b_{m,n}-2+i\}\\
 \Theta_2 \cap Y  = & \,  \{ \rho^{i}(s_0\theta(p,q)s_{p,q}) \, |\, a_{p,q}< a_{m,n}-3 +f(i) \mbox{ and }  b_{p,q} <   b_{m,n}-3+i\}
\end{align}
If one replaces every $<$  by  $\leq$ in the description above, one obtains the set 
\begin{equation*}
[\mbox{id}, \theta(m,n)s_{m,n}]\setminus \{\mbox{id}\} .
\end{equation*}
\end{rem}

% \begin{align}
%     X\cap [\mbox{id}, \theta(m,n) ]   = &\,  \{  \rho^i  (x_k)\, |\,  k\leq a_{m,n}+2-i \} \cup \{ \rho^i  (\sigma(x_k)) \, |\,  k\leq b_{m,n}+2-f(i)  \}. \\
%         \Theta \cap [\mbox{id}, \theta(m,n) ]  = &\,  \{\rho^{i}(\theta(p,q)) \, |\, a_{p,q}\leq a_{m,n}-i \mbox{ and }  b_{p,q} \leq b_{m,n}-f(i)\}.\\
%          \Theta_1 \cap [\mbox{id}, \theta(m,n) ]  = &\,   \{ \rho^{i}(\theta(p,q)s_{p,q})\, |\,  a_{p,q}< a_{m,n}-i \mbox{ and }  b_{p,q} < b_{m,n}-f(i)  \} \,  \cup  \\
%     & \,  \{ \rho^{i}(s_0\theta(p,q)) \, |\, a_{p,q}\leq a_{m,n}-3 +f(i) \mbox{ and }  b_{p,q} \leq   b_{m,n}-3+i\}.\\
%  \Theta_2 \cap [\mbox{id}, \theta(m,n) ]  = & \,  \{ \rho^{i}(s_0\theta(p,q)s_{p,q}) \, |\, a_{p,q}< a_{m,n}-3 +f(i) \mbox{ and }  b_{p,q} <   b_{m,n}-3+i\}.  
% \end{align}

\begin{figure}[ht]
  \begin{subfigure}[t]{.45\textwidth}
    \centering
    \includegraphics[trim=0cm 8cm 8cm 0cm,clip,width=\textwidth]{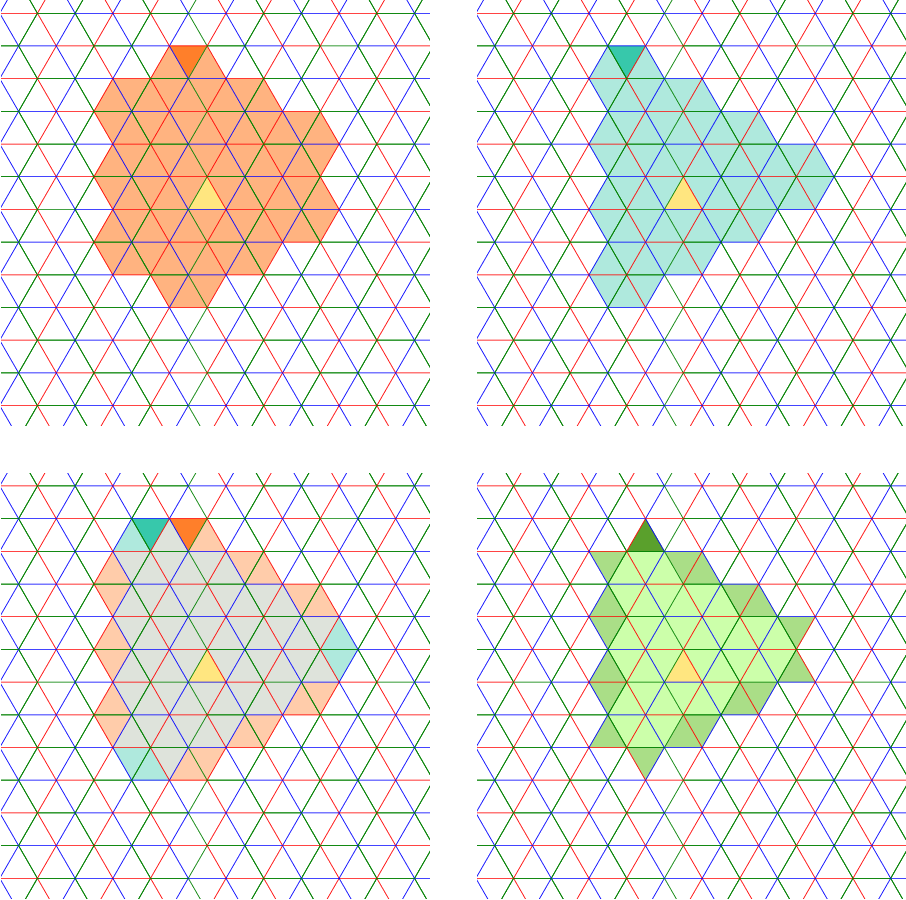}
    \caption{The interval $[\operatorname{id}, \theta(2,1)]$. The highlighted triangle corresponds to $\theta(2,1)$.}
    \label{fig: 3a}
  \end{subfigure}
  \quad
  \begin{subfigure}[t]{.45\textwidth} 
    \centering
    \includegraphics[trim=8cm 8cm 0cm 0cm,clip,width=\textwidth]{CICA2figs/LDLemainterseccion.eps}
    \caption{The interval $[\operatorname{id}, \theta(3,0)]$. The highlighted triangle corresponds to $\theta(3,0)$.}
    \label{fig: 3b}
  \end{subfigure}
    \smallbreak
  \begin{subfigure}[t]{.45\textwidth} 
    \centering
    \includegraphics[trim=0cm 0cm 8cm 8cm,clip,width=\textwidth]{CICA2figs/LDLemainterseccion.eps}
    \caption{The overlapping of Figure \ref{fig: 3a} and Figure \ref{fig: 3b}. The triangles in the intersection are colored with gray.}
    \label{fig: 3c}
  \end{subfigure}
  \quad
  \begin{subfigure}[t]{.45\textwidth} 
    \centering
    \includegraphics[trim=8cm 0cm 0cm 8cm,clip,width=\textwidth]{CICA2figs/LDLemainterseccion.eps}
    \caption{The interval $[\operatorname{id}, \theta(2,0)s_1]$. From lightest to darkest green:  The interval $[\operatorname{id}, \theta(2,0)]$, the set $\partial_{2,0}s_1\setminus\{\theta(2,0)s_1\}$ and the element $\theta(2,0)s_1$.}
    \label{fig: 3d}
  \end{subfigure}
  \caption{An example of Lemma \ref{lemma intersecccion} for the pair $(m,n)=(3,1)$.}
  \label{fig: interseccion}
\end{figure}

\begin{lem} \label{lemma intersecccion} 
Let $m$ and $n$ be positive integers.  Then, we have
\begin{equation} \label{eq interseccion intervalos}
    [\operatorname{id}, \theta(m-1,n)]   \cap  [\operatorname{id}, \theta(m,n-1)]  =  [\operatorname{id}, \theta(m-1,n-1)s_{m-1,n-1}].
\end{equation}
\end{lem}
\begin{proof}
%Let us consider the relative location of the triangles  $\theta(m-1,n)$, $\theta(m,n-1)$,  $\theta(m-1,n-1)$, and  $\theta(m-1,n-1)s$. These triangles are displayed in Figure \ref{fig: four triangles thetas}.  It follows that 
%\begin{equation}\label{eq interseccion convex}
%\mathcal{C}_{\theta(m-1,n)}\cap \mathcal{C}_{\theta(m,n-1)}=\mathcal{C}_{\theta(m-1,n-1)s}. 
%\end{equation}
This follows from Lemma \ref{menores que theta} and Lemma \ref{menores que s0theta}, see Figure \ref{fig: interseccion} for an illustration of this. 
\end{proof}

%\begin{equation*}
%\mathcal{C}_{\theta(m-1,n)}\cap \mathcal{C}_{\theta(m,n-1)}=[\operatorname{id}, \theta(m,n)] \uplus \partial_{m,n}s
%\end{equation*}

%\begin{figure}[h]
%    \centering
%    \centering
%    \includegraphics[scale=0.5]{CICA2figs/LDtres.eps}
%    \caption{From lightest to darkest gray: The relative location of the triangles corresponding to the elements $\theta(m-1,n-1) $, $\theta(m,n-1)$,\\ $\theta(m-1,n-1)s$, and $\theta(m-1,n)$.}
%    \label{fig: four triangles thetas}
%\end{figure}
% We know that  $W_f \lambda$ consists of $6$ points which form an hexagon with angles equal to $\pi/3$. The centroids of $\theta(m-1,n), \theta(m-1,n-1)s$, and $\theta(m-1,n)$ form an isosceles triangle with angles $2\pi/3, \pi/6$, and $\pi/6$. The same holds for the centroids of $w \theta(m-1,n), w \theta(m-1,n-1)s$, and $w \theta(m-1,n)$ for any $w\in W_f$.
%  Therefore
% \begin{equation*}
% \mathcal{C}_{\theta(m-1,n)}\cap \mathcal{C}_{\theta(m,n-1)}=\mathcal{C}_{\theta(m-1,n-1)s}.
% \end{equation*}
% By the geometric description in the Lemma \ref{menores que theta} and Lemma \ref{menores que s0theta} the statement follows.

% Let $\Phi^\vee$ be the coroot  data of $W$ with $\alpha^\vee, \beta^\vee$ the simple coroots. Note the centroid $\lambda=\operatorname{cen}(\theta(p,q))$ of the elements $\theta(p,q)$ satisfy
% \begin{align*}
% -\lfloor \langle\alpha^\vee, \lambda\rangle\rfloor&=p,\\
% -\lfloor \langle \beta^\vee, \lambda\rangle\rfloor&=q.
% \end{align*}

\subsection{The Hecke algebra of type \texorpdfstring{$\tilde{A}_2$}{A2 tilde}}\label{section Hecke}
The Hecke algebra $\mathcal{H}=\mathcal{H}(W)$ is the associative and unital $\mathbb{Z}[v, v^{-1}]$-algebra generated by  $\HH_{s_0}$, $\HH_{s_1}$ and  $\HH_{s_2}$ subject to the relations 
\begin{align*}
\HH_{s_i}^2 = (v^{-1}-v)\HH_{s_i}+1, \quad &\mbox{for } i\in \{0, 1, 2\},\\
\HH_{s_i}\HH_{s_j}\HH_{s_i} =\HH_{s_j}\HH_{s_i}\HH_{s_j},\quad &\mbox{for } i\neq j.
\end{align*}

Given $w\in W$ and a reduced expression $\underline{w}= r_1r_2\ldots r_k$ ($r_i\in S$, for all $i$) for $w$ we define $\HH_{w} =\HH_{r_1}\HH_{r_2}\cdots \HH_{r_k}  $. It is well-known that $\HH_{w}$ does not depend on the choice of the reduced expression.  The set $\{ \HH_{w} \, |\, w\in W  \}$ is a basis for $\mathcal{H}$ which is called the standard basis.  Let $\undH_{s} =\HH_{s}+v$. The right regular action of $\mathcal{H}$ is given by the formula

\begin{equation}  \label{multiply by b_s}
\HH_x\undH_s = \begin{cases}
\HH_{xs} +v\HH_{x}, & \text{if } xs>x;\\
\HH_{xs}+v^{-1}\HH_{x}, & \text{if } xs<x.
\end{cases}
\end{equation}

In their seminal paper \cite{KL79} Kazhdan and Lusztig introduced a new basis 
\begin{equation*}
\{\undH_{w} \, |\,  w\in W\}
\end{equation*} 
of $\mathcal{H}$ called the \emph{canonical basis}, also known as the \emph{Kazhdan-Lusztig basis}. We denote by $h_{x,w}(v)$ the corresponding Kazhdan-Lusztig polynomials defined by the equality 
\begin{equation}
    \undH_{w} = \sum_{x\in W} h_{x,w}(v) \HH_{x}.
\end{equation}
We stress that in this section we use Soergel's normalisation in \cite{soergel1997kazhdan} rather than the standard $q$-notation used by Kazhdan and Lusztig in \cite{KL79}.  The passage from Kazhdan-Lusztig polynomials $P_{x,w}(q)$ in their $q$-version to  Kazhdan-Lusztig polynomials $h_{x,w}(v)$ in their $v$-version is given by 
\begin{equation}\label{passage from v to q}
    h_{x,w}(v) = v^{\ell(w)-\ell (x) } P_{x,w}(v^{-2}).
\end{equation}
Henceforth, we use both versions of Kazhdan-Lusztig polynomials  without further notice.

For additional background material and an algorithm to compute the canonical basis, we refer the reader to \cite[Section 2]{soergel1997kazhdan}.  We conclude this section by recalling the formulas given in \cite{libedinsky2020affine}  for all the Kazhdan-Lusztig basis elements in $W$. 

For $x\in W$, we define 
\begin{equation*}
\NN_x\coloneqq \sum_{z\leq x}v^{\ell(x)-\ell(z)}\HH_{z}.
\end{equation*} 

\begin{prop}{\cite[Proposition 1.6]{libedinsky2020affine}}  \label{propo x}
%We have $\undH_{x_{n}}=\NN_{x_n}$ if $n\leq 3$ and $\undH_{x_4}= \NN_{x_4} +v\NN_{x_1}$. For $n\geq 5$ #
The following formula holds
\begin{equation}
   \label{KL x}    \undH_{x_n}  =     \begin{cases}
    \NN_{x_n}, & \mbox{ if $n\leq 3$;}\\
    \NN_{x_4}+v\NN_{x_{1}}, & \mbox{ if $n=4$;}\\
    \NN_{x_n}+v\NN_{x_{n-3}}, & \mbox{ if } n\geq 5 \mbox{ is odd;}\\
    \NN_{x_n}+v\NN_{x_{n-3}} +v\HH_{s_1s_0x_{n-5}}+v^2\HH_{s_0x_{n-5}} ,&   \mbox{ if } n\geq 5 \mbox{ is even.}\\
    \end{cases}  
\end{equation}
\end{prop}

\begin{prop}{\cite[Proposition 1.8]{libedinsky2020affine}} \label{formulas big region}
Let $m$ and $n$ be non-negative integers. We have
\begin{equation}
   \label{KL thetas}      \undH_{\theta(m,n)}  = \sum_{i=0}^{\min (m,n)} v^{2i}\NN_{\theta(m-i,n-i)}. 
\end{equation}
In particular, we have $\undH_{\theta(m,0)}=\NN_{\theta(m,0)}$ and $\undH_{\theta(0,n)}=\NN_{\theta(0,n)}$. Furthermore,  if $s=s_{m,n}$  we have
\begin{gather} \label{recursive formulas}
    \undH_{s_0}\undH_{\theta(m,n)} = \undH_{s_0\theta(m,n)}, \quad \undH_{\theta(m,n)}\undH_{s} = \undH_{\theta(m,n)s}, \\
\nonumber \mbox{ and } \quad \undH_{s_0}\undH_{\theta(m,n)}\undH_{s} =  \undH_{s_0\theta(m,n)s}.
\end{gather}
\end{prop}

Proposition \ref{propo x} and Proposition \ref{formulas big region} are enough to compute all the Kazhdan-Lusztig basis elements in view of the first claim in the following lemma. 

% In view of Lemma \ref{lemma LD is invariant}, it is enough to specify the formulas for the basis elements of the form $\undH_{x_n}$, $\undH_{\theta(m,n)}$, $\undH_{\theta(m,n)s_{m,n}}$, and $\undH_{s_0\theta(m,n)s_{m,n}}$. \\

\begin{lem} \label{lemma LD is invariant}
Suppose $\tau \in G$. Let $[x,y]$ be any Bruhat interval in $W$. Then, 
\begin{itemize}
  \item $P_{x,y}(q)=P_{\tau(x), \tau (y)}(q)$ (equivalently, $h_{x,y}(v)=h_{\tau(x), \tau (y)}(v)$).
  \item $[x,y]\simeq [\tau(x), \tau (y)] $ as posets.
\end{itemize}
\end{lem}
\begin{proof}
Both claims follow by the definition of Bruhat order and the definition of the Kazhdan-Lusztig polynomials.
\end{proof}

Lemma \ref{lemma LD is invariant} tells us that  if we prove the combinatorial invariance conjecture for a pair of intervals $[x,y]$ and $[x',y']$ then the conjecture is immediately true for any pair of intervals of the form  $[\tau (x),\tau(y)]$ and $[\tau'(x'),\tau'(y')]$, for all $\tau, \tau' \in G$. \\

\subsection{Trivial Kazhdan-Lusztig polynomials}\label{section trivialKL} 
For a Coxeter system $(W,S)$, define  $T=\bigcup_{w\in W} w S w^{-1}$ to be the set of all reflections. 
\begin{defn}
 Let $(W,S)$ be a Coxeter system. For $x\leq y$, we say that the interval $[x,y]$ satisfies  \emph{property P} if for every $z\in[x,y]$ we have
\begin{equation*}
\left|\{t\in T : z < tz\leq y\}\right|=\ell(y)-\ell(z).
\end{equation*}
\end{defn}

Thanks to the non-negativity theorem of the coefficients of the Kazhdan-Lusztig polynomials \cite{EWhodge} we can restate  \cite[Theorem C]{Carr91} as follows.

\begin{thm}\label{theorem C}
For a Coxeter system $(W,S)$, and $x\leq y$, the following are equivalent:
\begin{enumerate}
\item $P_{x,y}(q)=1$.
\item $P_{z,w}(q)=1$ for all $z\in [x,y]$.
\item The interval $[x,y]$ satisfies property P.
\end{enumerate}
\end{thm} 

\begin{defn}
The \emph{Bruhat graph} $\Omega_W([x,y])$ of $[x,y]$ is the directed graph defined as follows. The set of vertices is $[x,y]$ and the  set of edges is 
\begin{equation*}
\{(w, tw) : \mbox{ $w<tw$ and $t\in T$}\}.
\end{equation*}
\end{defn}

The following is \cite[Proposition 3.3]{Dy91}.

\begin{prop}\label{33dyer}
The Bruhat graph $\Omega_W([x,y])$ depends only on the poset type of the interval $[x,y]$.
\end{prop}

Putting all the results of this section together we have the following proposition.
\begin{prop}\label{prop trivialKL}
Let $(W,S)$ and $(W',S')$ be two Coxeter systems. Suppose $x,y\in W$ and $x',y'\in W'$ are elements such that $[x,y]\cong [x',y']$ as posets. If $P_{x,y}(q)=1$ then $P_{x',y'}(q)=1$.
\end{prop}
\begin{proof}
Suppose $P_{x,y}(q)=1$, by Theorem \ref{theorem C} we have $[x,y]$ satisfies the property P. By Proposition \ref{33dyer}, Property P is a poset invariant, therefore $[x',y']$ satisfies property P. By Theorem \ref{theorem C} we conclude $P_{x',y'}(q)=1$.
\end{proof}

\subsection{Monotonicity and content}\label{mono and content} First, we recall the monotonicity of Kazhdan-Lusztig polynomials proved in \cite{braden2001moment} for affine Weyl groups and in \cite{plaza2017graded} for arbitrary Coxeter systems. 

\begin{prop} \label{prop monotonicity}
If $x\leq z\leq y$ then 
\begin{equation}\label{v monotonicity}
h_{x,y}(v) - v^{\ell( z)-\ell (x) }   h_{z,y}(v) \in \bbN[v].  
\end{equation}
Equivalently,  
\begin{equation} \label{q monotonicity}
    P_{x,y}(q)-P_{z,y}(q) \in \bbN[q].
\end{equation}
\end{prop}

\begin{defn}
Given $h_1(v), h_2(v)\in \mathbb{Z}[v,v^{-1}]$, we write $h_1(v)\geq h_2(v)$ if $h_1(v)-h_2(v)\in \mathbb{N}[v,v^{-1}]$. Given $x\in W$ and $H\in \mathcal{H}$, we define $G_{x}(H)$ to be the coefficient of $\HH_{x}$ when we expand $H$ in terms of the standard basis of $\mathcal{H}$. We say an element $H\in \mathcal{H}$ is \emph{monotonic} if 
\begin{equation*}
G_{y}(H)\geq v^{\ell (x) -\ell (y)  } G_{x}(H),
\end{equation*}
for all $y\leq x$.  Finally, we write $H_1\leqH H_2$ if $G_x(H_{1} -H_{2}) \in \mathbb{N}[v,v^{-1}]$, for all $x\in W$. 
\end{defn}

An obvious example of a monotonic element is $\NN_{w}$. Also, by  \eqref{v monotonicity} we know that the canonical basis elements $\undH_{w}$ are monotonic in an arbitrary Coxeter system. On the other hand, it is easy to see that the sum of two monotonic elements is monotonic. The following (a weaker version of \cite[Lemma 3.3]{batistelli2021kazhdan}) is another way to produce new monotonic elements from old ones.

\begin{lem}\label{monotonic}
For $s\in S$, if $H\in \mathcal{H}$ is monotonic, then $H\undH_{s}$ is also monotonic. 
\end{lem}

\begin{defn}
Let $H\in \mathcal{H}$. We define the \emph{content} of $H$ as 
\begin{equation}
    c(H) = \sum_{x\in W} G_{x}(H)(1) \in \mathbb{Z}.
\end{equation}
\end{defn}

For instance, the content of $\NN_w$ is equal to the number of elements in the lower Bruhat interval 
$[\operatorname{id},w]$.  For $s=s_{m,n}$ we have the following formulas
\begin{align}
\label{less theta}    c(\NN_{\theta(m,n)}) & = \left|[\operatorname{id},\theta(m,n)]\right| = 3m^2+3n^2+12mn+9m+9n+6. \\
\label{less theta s}    c(\NN_{\theta(m,n)s})& = \left|[\operatorname{id},\theta(m,n)s]\right| = 
3m^2+3n^2+12mn+15m+15n+12. \\
\label{less s0 theta s}    c(\NN_{s_0\theta(m,n)s})& = \left|[\operatorname{id},s_0\theta(m,n)s]\right| = 
3m^2+3n^2+12mn+21m+21n+22. 
\end{align}
The formula in \eqref{less theta} is already  in \cite[Lemma 1.4]{libedinsky2020affine}. The remaining formulas can be  obtained in a similar way, using the geometric realization of $W$ as in Figure \ref{fig: four regions}. 
Using \eqref{multiply by b_s} one can easily prove the following lemma. 
\begin{lem} \label{double}
For any $H\in \mathcal{H}$  and $s\in S$ we have $c(H\undH_s)=2c(H). $   
\end{lem}

\begin{rem} \label{remark equality in Hecke}
We can use  the order $\leqH$ and the content to show an equality in $\mathcal{H}$. More precisely, let $H_1,H_2 \in \mathcal{H}$. If  $H_1\leqH  H_2 $ and  $c(H_1) = c(H_2) $ then $H_1=H_2$. 
\end{rem}

\section{New formulas for Kazhdan-Lusztig basis elements}  \label{section KL formulas}

Although Proposition \ref{formulas big region} is enough to compute all of the Kazhdan-Lusztig basis elements corresponding to elements in sets $\Theta_1$ and $\Theta_2$, we want more explicit formulas that allow us to compute the Kazhdan-Lusztig polynomials in a more direct way. This is the content of Proposition \ref{proposition KL thetas uno} and Proposition \ref{proposition KL theta dos} below. In this section we use the convention that $\undH_{\theta(m,n)}$ and $\NN_{\theta(m,n)}$ are zero if $m$ or $n$ are negative.

\begin{prop} \label{proposition KL thetas uno}
Let $m$ and $n$ be non-negative integers and $s=s_{m,n}$. We have
\begin{equation} \label{Kl thetas uno}
 \undH_{\theta(m,n)s}   = 
      \NN_{\theta(m,n)s}  +   v\undH_{\theta(m-1,n)} + v \undH_{\theta(m,n-1).}
\end{equation}

\end{prop}
\begin{proof}
The  $m=n=0$ case is easily checked by hand. \\
Suppose that $m>0$ and $n=0$ (the case $m=0$ and $n>0$ is  similar). By Proposition \ref{formulas big region} we have $\undH_{\theta(m,0)}=\NN_{\theta(m,0)}$. Therefore, Equation \eqref{recursive formulas} implies $\undH_{\theta(m,0)s}=\NN_{\theta(m,0)}\undH_{s}$. Therefore, this case is reduced to check the identity
\begin{equation}
    \NN_{\theta(m,0)}\undH_{s} = \NN_{\theta(m,0)s} +v \NN_{\theta(m-1,0)}. 
\end{equation}
To prove this we use Remark \ref{remark equality in Hecke}. Let 
\begin{equation*}
H_1= \NN_{\theta(m,0)}\undH_{s} \mbox{ and } H_2=\NN_{\theta(m,0)s} +v \NN_{\theta(m-1,0)}.
\end{equation*}
By Lemma \ref{double}, Equation \eqref{less theta} and Equation \eqref{less theta s} we obtain 
\begin{equation}
    c(H_1 ) = c( H_2) = 2(3m^2+9m+6). 
\end{equation}

We have
\begin{equation} \label{eq H2 coeff}
    G_{x}(H_2) = \begin{cases}
    v^{\ell (\theta(m,0)s)-\ell(x)},  & \mbox{if } x \leq \theta(m,0)s \mbox{ and } x\not \leq \theta(m-1,0);\\
    v^{\ell (\theta(m,0)s)-\ell(x)}+ v^{\ell (\theta(m-1,0))+1 -\ell(x)}, & \mbox{if } x \leq \theta(m-1,0);\\
    0, &\mbox{otherwise.}
  \end{cases}
\end{equation}
On the other hand, by Equation \eqref{multiply by b_s} we get
\begin{equation} \label{eq lala}
 G_{\theta(m,0)s}(H_1) = 1 \quad \mbox{ and } \quad 
G_{\theta(m-1,0)}(H_1)=v^3+v.  
\end{equation}
By Lemma \ref{monotonic} $H_1$ is monotonic. Therefore,  \eqref{eq lala} implies 
\begin{equation}\label{eq H1 coeff}
    G_{x}(H_1) \geq  \begin{cases}
    v^{\ell (\theta(m,0)s)-\ell(x)},  & \mbox{if } x \leq \theta(m,0)s \mbox{ and } x\not \leq \theta(m-1,0);\\
    v^{\ell (\theta(m,0)s)-\ell(x)}+ v^{\ell (\theta(m-1,0))+1 -\ell(x)}, & \mbox{if } x \leq \theta(m-1,0).
  \end{cases}
\end{equation}
In the lower part of \eqref{eq H1 coeff} we are using that $\ell( \theta(m-1,0) )  +3   =   \ell( \theta(m,0)s )$. By combining \eqref{eq H2 coeff} and \eqref{eq H1 coeff} we obtain  $H_1\leqH H_2$, thus proving the lemma in this case.\\

We now assume that $m,n>0$. By \eqref{KL thetas} we have
\begin{equation} \label{decomposition Theta equal N plus Theta}
    \undH_{\theta(m,n)} = \NN_{\theta(m,n) } + v^2 \undH_{\theta(m-1,n-1)}.
\end{equation}
Multiplying on the right by $\undH_s$ (note that $s=s_{m,n}=s_{m-1,n-1}$) and assuming by induction that \eqref{Kl thetas uno}  holds for $\theta(m-1,n-1)$ we obtain that
\begin{align*}
     \undH_{\theta(m,n)s } &  = \NN_{\theta(m,n) }\undH_{s} + v^2 \undH_{\theta(m-1,n-1)}\undH_{s}  \\
    &  = \NN_{\theta(m,n) }\undH_{s} + v^2 \undH_{\theta(m-1,n-1)s}  \\ 
     &  = \NN_{\theta(m,n) }\undH_{s} + v^2\NN_{ \theta(m-1,n-1)s }  +v^3 \undH_{ \theta(m-2,n-1) }+v^3 \undH_{ \theta(m-1,n-2) }.
\end{align*}

Therefore, using \eqref{decomposition Theta equal N plus Theta} twice, our claim reduces to prove the following identity
\begin{equation} \label{lo que falta}
    \NN_{\theta(m,n)}\undH_{s} +v^2\NN_{\theta(m-1,n-1)s} = \NN_{\theta(m,n)s} +v\NN_{\theta(m-1,n)} +v\NN_{\theta(m,n-1)}. 
\end{equation}
As before, we use Remark \ref{remark equality in Hecke}. Let $L$ and $R$ be the left-hand side and the right-hand side of \eqref{lo que falta}, respectively. A combination of \eqref{less theta}, \eqref{less theta s}, and Lemma \ref{double} yields 
\begin{equation}  \label{the contents are equal}
    c(L)=c(R)= 3(3m^2+3n^2+12mn+5m+5n+4). 
\end{equation}
It remains to show that $L\leqH R$. Set $a_0=\theta(m,n)s$, $a_1=\theta(m-1,n)$, $a_2=\theta(m,n-1)$, and $a_3=\theta(m-1,n-1)s$. By the definition of $R$ and  Lemma \ref{lemma intersecccion} we get

\begin{equation} \label{eq R coeff}
    G_{x}(R) = \begin{cases}
    v^{\ell (a_0)-\ell(x)},  & \mbox{if } x \leq a_0, \, x\not \leq a_1, \mbox{ and }  x\not \leq a_2 ;\\
    v^{\ell (a_0)-\ell(x)}+ v^{\ell (a_1)+1 -\ell(x)} , & \mbox{if } x \leq a_1 \mbox{ and } x \not \leq a_2 ;\\
   v^{\ell (a_0)-\ell(x)}+ v^{\ell (a_2)+1 -\ell(x)} , & \mbox{if } x \not \leq a_1 \mbox{ and } x  \leq a_2 ;\\
      v^{\ell (a_0)-\ell(x)}+ 2v^{\ell (a_2)+1 -\ell(x)} , & \mbox{if } x  \leq a_3  ;\\
    0, &\mbox{otherwise.}
  \end{cases}
\end{equation}
We have used that $\ell(a_1)=\ell(a_2)$. By the definition of $L$ and  Equation \eqref{multiply by b_s} we get
 \begin{align}
     G_{a_0}(L) & =1;\\
     G_{a_1}(L)&= v^3+v; \\ 
     G_{a_2}(L)&= v^3+v; \\ 
     G_{a_3}(L) &= v^4+2v^2.  
 \end{align}
Using Lemma \ref{monotonic} and the fact that the sum of two monotonic elements is monotonic we conclude that $L$ is monotonic as well. It follows that 
\begin{equation} \label{eq L coeff}
    G_{x}(L) \geq  \begin{cases}
    v^{\ell (a_0)-\ell(x)},  & \mbox{if } x \leq a_0, \, x\not \leq a_1, \mbox{ and }  x\not \leq a_2 ;\\
    v^{\ell (a_0)-\ell(x)}+ v^{\ell (a_1)+1 -\ell(x)} , & \mbox{if } x \leq a_1 \mbox{ and } x \not \leq a_2 ;\\
   v^{\ell (a_0)-\ell(x)}+ v^{\ell (a_2)+1 -\ell(x)} , & \mbox{if } x \not \leq a_1 \mbox{ and } x  \leq a_2 ;\\
      v^{\ell (a_0)-\ell(x)}+ 2v^{\ell (a_2)+1 -\ell(x)} , & \mbox{if } x  \leq a_3 ,
  \end{cases}
\end{equation}
Where we have used the following equalities.
\begin{align}
    \ell( a_1 )  +3  & =   \ell( a_0 ), \\
      \ell( a_2 ) +3 & =   \ell( a_0 ),\\
        \ell( a_3 ) +4 & =   \ell( a_0 ).
\end{align}
  By combining \eqref{eq R coeff} and \eqref{eq L coeff} we obtain  $L\leqH R$, as we wanted to show.
\end{proof}

The formula for $ \undH_{s_0\theta(m,n)s_{m,n}}$ is  more involved. %Furthermore, it admits certain symmetry that will be important for us in the forthcoming section. 
To describe this formula we need to consider the following element. For $x,y\in W$,  we define
\begin{equation}
   \MM_{x,y}\coloneqq \sum_{ w \leq x \mbox{ or } w\leq y} v^{\ell (x)-\ell (w)} \HH_{w}.   
\end{equation}
\begin{rem}\label{Mxy=Myx}
One has the equality $\MM_{x,y}= \MM_{y,x}$ if and only if $\ell (x)=\ell (y)$. 
\end{rem}

The content of $\MM_{x,y}$ equals the number of elements in the union of  the lower intervals $[\operatorname{id},x]$ and $[\operatorname{id},y]$. In particular, for any pair $(m,n)$ of positive integers we have 
\begin{equation}  \label{content M }
   c(\MM_{\rho (\theta(m,n-1))s, \rho^2 (\theta(m-1,n) )s })   = 3m^2+3n^2+12mn+9m+9n+2,
\end{equation}
where $s=s_{m,n}$. Note that $s$ is such that $\rho (\theta(m,n-1))s>\rho (\theta(m,n-1))$ and $\rho^2 (\theta(m-1,n) )s>\rho^2 (\theta(m-1,n) )$.

\begin{prop}  \label{proposition KL theta dos}
Let $m$ and $n$ be non-negative integers and $s=s_{m,n}$. \\\\
 If $m=n=0$ then 
 \begin{equation} \label{eq theta dos 00}
   \undH_{s_0\theta(0,0)s} = \NN_{s_0\theta(0,0)s} + v^2\NN_{s_0}.   
 \end{equation}
 If $m>0 $ and $n=0$ then 
      \begin{align}
            \label{caso m>0 n=0 verison 1}     \undH_{s_0\theta(m,0)s} & = \NN_{s_0\theta(m,0)s} + v\MM_{s_0\theta(m-1,0),\rho (\theta(m-1,0) ) }  + v \undH_{ \rho^2( \theta(m-1,0))s}, \\
             \label{caso m>0 n=0 verison 2}                            & = \NN_{s_0\theta(m,0)s}  + v\MM_{\rho^2( \theta(m-1,0))s, \rho (\theta(m-1,0))}  + v \undH_{ s_0\theta(m-1,0)}.     
      \end{align}
 If $m=0 $ and $n>0$ then 
    \begin{align}
   \label{caso m=0 n>0 verison 1}    \undH_{s_0\theta(0,n)s} & = \NN_{s_0\theta(0,n)s} + v\MM_{s_0\theta(0,n-1),\rho^2 (\theta(0,n-1) ) }  + v \undH_{ \rho( \theta(0,n-1))s}, \\
         \label{caso m=0 n>0 verison 2}                            & = \NN_{s_0\theta(0,n)s}  + v\MM_{\rho( \theta(0,n-1))s, \rho^2 (\theta(0,n-1))}  + v \undH_{ s_0\theta(0,n-1)}.
    \end{align}  
If $m>0 $ and $n>0$ then 
     \begin{align}
    \label{caso m>0 n>0 verison 1}     \undH_{s_0\theta(m,n)s} &  = \NN_{s_0\theta(m,n)s} +  
        v\MM_{s_0\theta(m,n-1), s_0\theta(m-1,n)} +v \undH_{\rho(\theta(m,n-1))s} + v \undH_{\rho^2(\theta(m-1,n))s}, \\
    \label{caso m>0 n>0 verison 2}                                &  = \NN_{s_0\theta(m,n)s}  +v\MM_{\rho(\theta(m,n-1))s,\rho^2(\theta(m-1,n))s  } +v \undH_{  s_0\theta(m,n-1)} + v \undH_{s_0\theta(m-1,n)  }. 
     \end{align}   
\end{prop}
\begin{proof}
We only prove the case $m>0$ and $n>0$, the proof of the remaining cases being analogous. 

Multiplying \eqref{Kl thetas uno} by $\undH_{s_0}$ on the left and using \eqref{recursive formulas} we obtain
% By Proposition \ref{proposition KL thetas uno} we know that
% \begin{equation}
%  \undH_{\theta(m,n)s} = \NN_{\theta(m,n)s} +v \undH_{\theta(m-1,n)} + v\undH_{\theta(m,n-1)}.   
% \end{equation}
% Multiplying this equality by $\undH_{s_0}$ on the left and using \eqref{recursive formulas} we obtain
\begin{equation}
 \undH_{s_0\theta(m,n)s} =\undH_{s_0} \NN_{\theta(m,n)s} +v \undH_{s_0\theta(m-1,n)} + v\undH_{s_0\theta(m,n-1)}.   
\end{equation}
Therefore, the proof of \eqref{caso m>0 n>0 verison 2} reduces to show the identity
\begin{equation}\label{loquefalta2}
   \undH_{s_0} \NN_{\theta(m,n)s} = \NN_{s_0\theta(m,n)s} +v\MM_{\rho(\theta(m,n-1))s,\rho^2(\theta(m-1,n))s}.
\end{equation}
To prove this we use Remark \ref{remark equality in Hecke}.  Let $L$ and $R$ be the left-hand side and the right-hand side of \eqref{loquefalta2}, respectively. By \eqref{less theta s}, \eqref{less s0 theta s}, \eqref{content M }, and (the left version of) Lemma \ref{double} we get
\begin{equation*}
c(L)=c(R)=6m^2+6n^2+24mn+30m+30n+24.
\end{equation*}
It remains to show that $L\leqH R$. Set $a_0=s_0\theta(m,n)s$, $a_1=\rho(\theta(m,n-1))s$, and $a_2=\rho^2(\theta(m-1,n))s$. By the definition of $R$ we get
\begin{equation} \label{GxR}
    G_{x}(R) = \begin{cases}
    v^{\ell (a_0)-\ell(x)},  & \mbox{if } x\leq a_0, \, x \not\leq a_1 \mbox{ and } x \not\leq a_2;\\
    v^{\ell (a_0)-\ell(x)}+v^{\ell (a_1)-\ell(x)+1},  & \mbox{if } x \leq a_1 \mbox{ or } x \leq a_2;\\
    0, &\mbox{otherwise.}
  \end{cases}
\end{equation}
In the second line, we have used that $\ell(a_1)=\ell(a_2)$. We notice that $s_0 a_1<a_1$ and $s_0 a_2< a_2$. Then, the left version of Equation \eqref{multiply by b_s} yields $G_{a_0}(L)= 1$ and $G_{a_1}(L)=G_{a_2}(L)= v^3+v$. 

By Lemma \ref{monotonic} $L$ is monotonic. It follows that 
\begin{equation} \label{GxL}
    G_{x}(L) \geq  \begin{cases}
    v^{\ell (a_0)-\ell(x)},  & \mbox{if } x\leq a_0, \, x \not\leq a_1, \mbox{ and } x \not\leq a_2;\\
    v^{\ell (a_0)-\ell(x)}+v^{\ell (a_1)-\ell(x)+1},  & \mbox{if } x \leq a_1 \mbox{ or } x \leq a_2.
  \end{cases}
\end{equation}
Where we have used that $G_{a_1}(L)=G_{a_2}(L)$, and the fact that $\ell(a_0)=\ell(a_1)+3$. By combining \eqref{GxR} and  \eqref{GxL} we obtain  $L\leqH R$. This finishes the proof of \eqref{caso m>0 n>0 verison 2}.\\
% Note that $a_1>s_0 a_1$, to be more precise $\rho^2(s_0\theta(m,n-2))s$ gives a reduced expression of $s_0a_1$ if $n\geq 2$ and $\rho(x_{2m-2})$ gives a reduced expression of $s_0a_1$ if $n=1$ (we set the convention $x_0=\operatorname{id}$ here). Similarly $a_2>s_0 a_2$, to be more precise $\rho(s_0\theta(m-2,n))s$ gives a reduced expression of $s_0a_2$ if $m\geq 2$ and $\rho^2(x_{2n-2})$ gives a reduced expression of $s_0a_2$ if $m=1$ (with the same convention for $x_0$ as before). By Equation \eqref{multiply by b_s} and the previous two facts we obtain
%  \begin{align}
%      G_{a_0}(L)&= 1;\\
%      G_{a_1}(L)&= v^3+v; \\ 
%      G_{a_2}(L)&= v^3+v.
%  \end{align}
%%Gaston is writing above

Finally, \eqref{caso m>0 n>0 verison 1} is obtained from \eqref{caso m>0 n>0 verison 2} by inverting and then applying a power of $\rho$. Indeed,  we have the following identity in $W$:
\begin{equation}
    \left[  \rho^j(\theta(m,n))  \right]^{-1} = \rho^{j+n-m}(\theta(n,m)). 
\end{equation}
Using this identity to invert all the elements of $W$ occurring in \eqref{caso m>0 n>0 verison 2} we get
\begin{align}
     \undH_{s\rho^{k} (\theta(n,m)) s_0}  = &\NN_{s\rho^{k} (\theta(n,m)) s_0} +v\MM_{ s\rho^{k} (\theta(n-1,m)) ,s \rho^{k}(\theta(n,m-1) ) } +\\
     \nonumber & v\undH_{\rho^{k-1} (\theta(n-1,m) ) s_0 } + v\undH_{ \rho^{k+1} (\theta(n,m-1)) s_0 }, 
\end{align}
where $k=n-m$. Then, we act by $\rho^{-k}$ on the equality above, and using the fact that $\rho^{-k}(s)=s_0$ and $\rho^{-k}(s_{0})=s_{n,m} \eqqcolon \tilde{s} $, we obtain
\begin{equation*}
   \undH_{s_0\theta(n,m) \tilde{s}}  =  \NN_{s_0\theta(n,m) \tilde{s}}+v\MM_{s_0\theta(n,m-1),s_0\theta(n-1,m)} + v \undH_{\rho^2(\theta(n-1,m))\tilde{s}} + v\undH_{\rho (\theta(n,m-1)) \tilde{s}}.     
\end{equation*}
Since $\ell(s_0\theta(n,m-1))=\ell(s_0\theta(n-1,m))$, by Remark \ref{Mxy=Myx} we get
\begin{equation*}
   \undH_{s_0\theta(n,m) \tilde{s}}  =  \NN_{s_0\theta(n,m) \tilde{s}}+v\MM_{s_0\theta(n-1,m),s_0\theta(n,m-1)} + v \undH_{\rho^2(\theta(n-1,m))\tilde{s}} + v\undH_{\rho (\theta(n,m-1)) \tilde{s}}.     
\end{equation*}
This last equality is  \eqref{caso m>0 n>0 verison 1}  with the roles of $m$ and $n$ switched. 
\end{proof}

\section{Two poset invariants}\label{inv}

For the proof of the main result of this paper it will be useful to have at hand certain invariants that will allow us to easily discard the occurrence of certain poset isomorphisms. 
\subsection{\texorpdfstring{$m$}{m}-joins}
%One of these invariants is the set of $m$-joins of a pair of elements in a given interval. 
\begin{defn}
Let $[x,y]$ be an interval and $m$ be a positive integer. Let $a,b\in [x,y]$ be such that $\ell (a) =\ell (b)$. We say that an element $z\in[x,y]$ is an \emph{$m$-join} of $a$ and $b$ if $a\leq z$, $b\leq z$, and $\ell(z) =\ell (a)+m=\ell (b)+m$. We denote by $\joins{a}{b}{x}{y}{m}$  the set of all $m$-joins of $a$ and $b$ in the interval $[x,y]$.
\end{defn}

\begin{rem}\label{1joins}
It is worth mentioning  the fact that $1$-joins have already been relevant for the combinatorial invariance conjecture. As a matter of fact, it was proved in  \cite[Theorem 3.2]{BCM06} that for $a\neq b$ we have  
\begin{equation} \label{1 joins equation}
\left| \joins{a}{b}{x}{y}{1} \right|\leq 2
\end{equation}
 in any Coxeter system. Although we will not use this result in our proof, it is interesting to notice the similarity between  \eqref{1 joins equation} and  our Lemma \ref{Joins of theta two}.
\end{rem}

The following result is immediate from the definitions. 

\begin{lem}\label{lemma Join to Join}
Let  $\phi\colon [x,y]\longrightarrow [x',y']$ be an isomorphism of posets. Suppose $a,b\in [x,y]$ are such that $\ell (a) = \ell (b)$. Then,  $\phi (\joins{a}{b}{x}{y}{m} ) = \joins{\phi(a)}{\phi(b)}{x'}{y'}{m}$, for all $m$. In particular, $\joins{a}{b}{x}{y}{m}$ and $\joins{\phi(a)}{\phi(b)}{x'}{y'}{m}$ have the same number of elements.
\end{lem}

As we have seen in Proposition \ref{proposition KL theta dos}, there are certain key elements in $W$ that allow us to compute the canonical basis elements in a simple way, namely those indexing the $\MM$ and $\undH$ elements appearing in the formulas in Proposition \ref{proposition KL theta dos}. In the following lemma, we will explore the number of $2$-joins of some  of these elements. These numbers are going to be important poset invariants to be used in our proof of the conjecture.

\begin{figure}[ht]
    \centering
\includegraphics[width=0.35\textwidth]{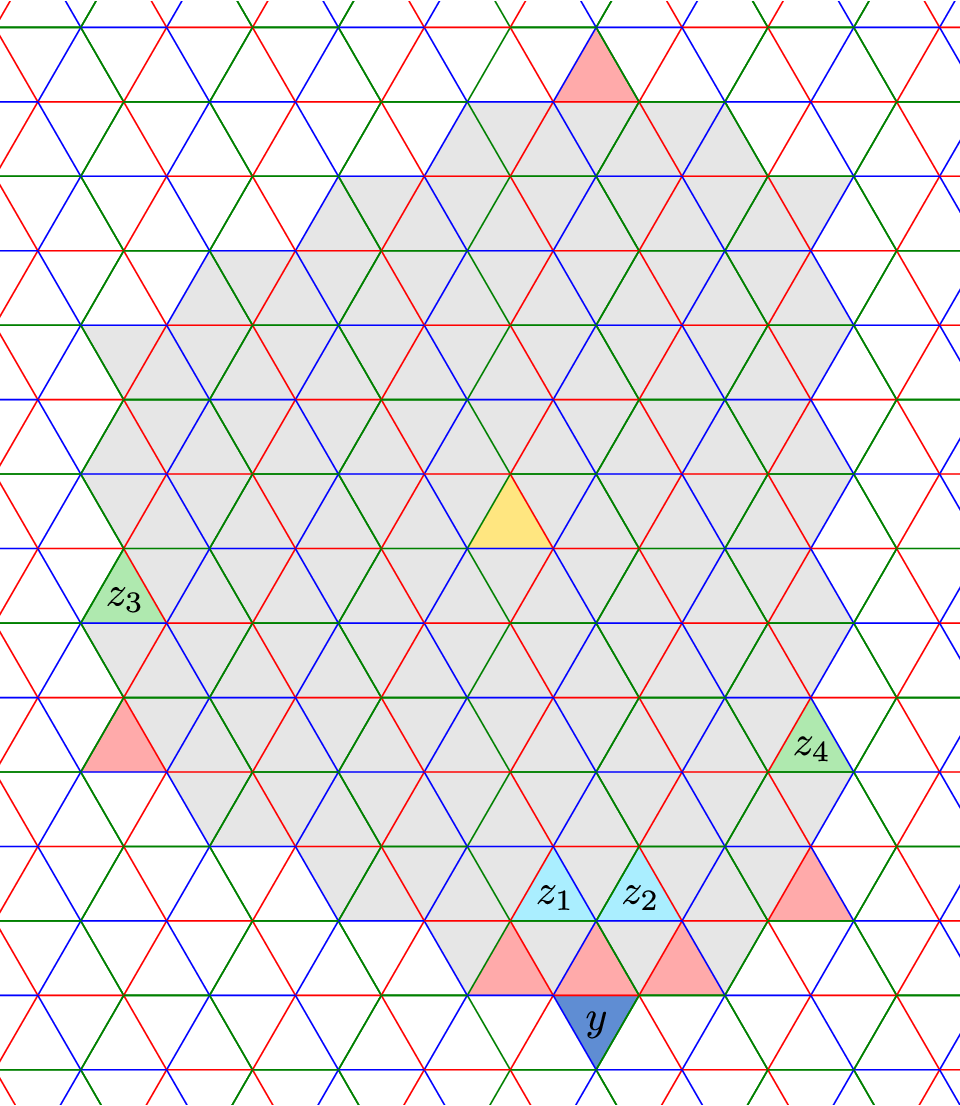}\\
    \caption{The triangles corresponding to the elements of the set $C_y$ for $y=s_0\theta(1,3)s_2$ are colored in pink. For reference, we colored all alcoves inside the interval $[\operatorname{id},y]$. As before, the yellow triangle corresponds to the identity in $W$.}
    \label{fig: lema28}
\end{figure}

\begin{lem}\label{Joins of theta two}
Let $[x,y]$ be an interval such that $y=s_0\theta(m,n)s, $ where $s=s_{m,n}$. When they make sense, we define
\begin{align} \label{zetas}
z_1&\coloneqq s_0\theta(m,n-1); 
&z_2&\coloneqq s_0\theta(m-1,n); \\
z_3&\coloneqq \rho^2(\theta(m-1,n))s; 
&z_4&\coloneqq \rho(\theta(m,n-1))s.
\end{align}
We remark that $z_1$ and $z_4$ (resp. $z_2$ and $z_3$) are only defined if $n>0$ (resp. $m>0$). %If $n=0$ (resp. $m=0$) then the elements $z_1$ and $z_4$ (resp. $z_2$ and $z_3$) are neglected.
Suppose $i$ and $j$ are such that both $z_i$ and $z_j$ are defined and belong to $[x,y]$ with $i\neq j$. We have
 \begin{equation*}
\left| \joins{z_i}{z_j}{x}{y}{2} \right|=
\begin{cases} 
3, & \mbox{ if $\{i,j\}=\{1,2\}$ or $\{i,j\}=\{3,4\}$;} \\
2, & \mbox{ otherwise}.
\end{cases}
\end{equation*}
\end{lem}
\begin{proof} 
We only prove the case where $m>0$ and $n>0$. This is when the four elements $z_1$, $z_2$, $z_3$, and $z_4$ are defined. The remaining cases ($m=0 $ and $n>0$, or  $m>0$ and $n=0$) are similar and easier.

\noindent
We notice that $\ell(z_i)=\ell(y)-3$ for $i\in\{1,2,3,4\}$. Therefore, 
\begin{equation} \label{joins written in coatoms}
  \joins{z_i}{z_j}{x}{y}{2} =\{ w\in C_y \, |\,   z_i< w \mbox{ and } z_j < w\}.  
\end{equation}
We recall that the set $C_y$ is given in Lemma \ref{corolario ganador}. See Figure \ref{fig: lema28} for a picture of the elements of $C_y$ and the $z_i$'s. The elements of $C_y$ relate with the $z_i$'s in the following way:
\begin{align}
z_3, z_4      &< \theta(m,n)s;\\
z_2, z_3, z_4 &< \rho (\theta(m-1,n+1))s;\\
z_1, z_2, z_3 &< s_0\theta(m-1,n+1);\\  \label{zetas and Cy}
z_1, z_3, z_4 &< \rho^2 (\theta(m+1,n-1))s;\\ 
z_1, z_2, z_4 &< s_0\theta(m+1,n-1);\\
z_1, z_2      &< s_0\theta(m,n).
\end{align}
The lemma follows by combining \eqref{joins written in coatoms} and \eqref{zetas and Cy}. For instance, 
\begin{equation}
    \joins{z_1}{z_2}{x}{y}{2}= \{s_0\theta(m-1,n+1) ,  s_0\theta(m+1,n-1) ,s_0\theta(m,n)  \}
\end{equation}
since $s_0\theta(m-1,n+1)$, $ s_0\theta(m+1,n-1)$ and $s_0\theta(m,n)$  are the only $3$ elements of $C_y$ which are simultaneously greater than  $z_1$ and $z_2$.
\end{proof}

\subsection{Z-invariants}

  The following result \cite[Exercises 5.7 and 5.8]{BB05} is fundamental for our proof, as explained in the introduction.

\begin{prop}\label{proposition menor que cuatro}
If $[x,y]\sim [x',y']$ and $\ell(y)-\ell(x)\leq 4$ then $P_{x,y}(q)=P_{x',y'}(q)$. 
\end{prop}

Before embarking in the proof, we need to define the key poset invariant.

\begin{defn}
Let $[x,y] $ be any interval. For any integer $m\geq 3$ we define
\begin{equation}
     Z_{x,y}^m = \{ z\in [x,y] \, |\, \ell(y) -\ell(z)=m \mbox{ and } P_{z,y}(q)=1+q  \}.
\end{equation}
\end{defn}

\begin{lem} \label{lemma Z to Z}
Let  $\phi\colon [x,y]\longrightarrow [x',y']$ be a poset isomorphism. If $1\leq m\leq 4$,  then  $\phi(Z_{x,y}^m) =Z_{x',y'}^m$. In particular, $|Z_{x,y}^m| =|Z_{x',y'}^m|$.
\end{lem}
\begin{proof}
This is a direct consequence of Proposition \ref{proposition menor que cuatro}. %together with the Chain Property \cite[2.6]{MR644668}. 
\end{proof}

The following lemmas show the power of the Z-invariant. 

\begin{lem}\label{wrong lemma uno}
Let $[x,y]$ be an interval with $y\in \Theta_1 \cup\Theta_2 \cup X$. If $|Z_{x,y}^3|=1$ then $P_{x,y}(q)=1+q$. 
\end{lem}
\begin{proof}
 We need to split the proof in three cases: $y\in\Theta_1$, $y\in \Theta_2$, or $y\in X$.\\\\
Let us first assume that $y\in \Theta_1$. We can assume that $y=\theta(m,n)s_{m,n}$ for some $m,n\in \mathbb{N}$. By Proposition \ref{proposition KL thetas uno} we have $Z_{x,y}^3=\{\theta(m-1,n) \}$ or $Z_{x,y}^3=\{\theta(m,n-1) \}$. Suppose we are in the former case (the latter is treated similarly and is going to be omitted). Then $m>0$ and  $\theta(m-1,n)\in [x,y]$, and $\theta(m,n-1) \notin [x,y] $ (if $\theta(m,n-1)\in [x,y]$ then we would have $\theta(m,n-1)\in Z^{3}_{x,y}$, and therefore $|Z^{3}_{x,y}|=2$, contradicting our hypothesis).\\
Suppose that $P_{x,y}(q)\neq 1+q$. An inspection of \eqref{KL thetas} and \eqref{Kl thetas uno} reveals that $m>1$, $n>0$ and $\theta(m-2,n-1)\in [x,y]$. In particular, we have $x \leq  \theta(m-2,n-1)  $. Since $\theta(m-2,n-1) < \theta(m,n-1)$ we conclude that $\theta(m,n-1)\in [x,y] $ contradicting our conclusion in the last paragraph. Thus, $P_{x,y}(q)=1+q$, proving the lemma in this case. \\\\
Suppose now that $y\in \Theta_2 $. We can assume that $y=s_0\theta(m,n)s_{m,n}$.  The result follows by a case-by-case inspection of the formulas given in Proposition \ref{proposition KL theta dos}. Indeed, if $m=n=0$ then $Z_{x,y}^3=\emptyset$ and there is nothing to prove.
Suppose now that $m>0$ and $n=0$. By \eqref{caso m>0 n=0 verison 1} or \eqref{caso m>0 n=0 verison 2}  we have   $Z_{x,y}^3=\{ s_0\theta(m-1,0) \} $ or $Z_{x,y}^3=\{ \rho^2( \theta(m-1,0))s_{m,n}\}  $. If $Z_{x,y}^3=\{ s_0\theta(m-1,0) \} $ (resp. $Z_{x,y}^3=\{ \rho^2( \theta(m-1,0))s_{m,n}\}  $) then we use \eqref{caso m>0 n=0 verison 1} (resp. \eqref{caso m>0 n=0 verison 2}) to conclude that $P_{x,y}(q)=1+q$. The case $m=0$ and $n>0$ is treated similarly.\\
We can now assume that $m>0$ and $n>0$. By \eqref{caso m>0 n>0 verison 1} or \eqref{caso m>0 n>0 verison 2} we have that $Z_{x,y}^{3}=\{z_i\}$, for some $1\leq i\leq 4$ and where $z_i$ is as in Lemma \ref{Joins of theta two}. If $Z_{x,y}^3=\{z_1  \}$ or $Z_{x,y}^3=\{z_2 \}$  (resp. $Z_{x,y}^3=\{z_3  \}$ or $Z_{x,y}^3=\{z_4 \}$) then we use \eqref{caso m>0 n>0 verison 1}  (resp.  \eqref{caso m>0 n>0 verison 2}) to conclude that $P_{x,y}(q)=1+q$, as we wanted to show.\\\\
Finally, we suppose $y\in X$. Without loss of generality, $y=x_k$. Furthermore, we assume that $k$ is even and greater than or equal to $6$, since the other cases are easier. We prove something  stronger, namely that if $Z_{x,y}^3\neq \emptyset$ then $P_{x,y}(q)=1+q$. To see this we first notice that $x_{k-3} $ and $s_1s_0x_{k-5}$ (resp. $s_0x_{k-5}$) are incomparable in the Bruhat order. Furthermore, a  direct computation reveals that if $x< s_1s_0x_{k-5}$ and $x\neq s_0x_{k-5},$ then $x\leq x_{k-3}$. On the other hand, by Equation \eqref{KL x}  if $Z_{x,y}^3\neq \emptyset$ then we have that either $x\leq x_{k-3}$ or $x\leq s_1s_0x_{k-5}$.  We must have one and only one of the following possibilities: $x\leq x_{k-3}$, $x=s_1s_0x_{k-5}$ or $x=s_0x_{k-5}$. An inspection of the formulas in Equation \eqref{KL x} shows that in any of these three cases we have $P_{x,y}(q)=1+q$.
\end{proof}

\begin{rem}  \label{remark 1 y 1+q en X}
 The argument given in the final case of the proof of Lemma \ref{wrong lemma uno} shows that if  $ y\in X$ and $x\leq y$ we have
 \begin{equation}
     P_{x,y}(q)= \begin{cases}
     1, & \mbox{if } Z_{x,y}^3=\emptyset;\\
     1+q,& \mbox{if } Z_{x,y}^3\neq \emptyset.
     \end{cases}
 \end{equation}
\end{rem}

\begin{lem}\label{lemma Theta 2 with 1 in Z}
Let $[x,y]$ be an interval with $y\in \Theta_1  \cup X$. If $Z^3_{x,y}=\emptyset$ then $P_{x,y}(q)=1$. 
\end{lem}
\begin{proof}
The claim  follows by  an inspection of the formulas in Proposition \ref{propo x} and Proposition \ref{proposition KL thetas uno}. 
\end{proof}

\begin{lem}\label{wrong lemma dos}
Let $[x,y]$ be an interval with $y\in \Theta_2$. Suppose that $Z_{x,y}^3=\emptyset$. If $P_{x,y}(q)\neq 1$ then the following three conditions hold:
\begin{enumerate}
    \item $P_{x,y}(q)=1+q$;
    \item $|Z_{x,y}^4|=1$;
    \item $\ell(y)-\ell(x) =4$ or $\ell(y)-\ell(x)=5$. 
\end{enumerate}
\end{lem}

\begin{proof}
We can assume that $y= s_0 \theta(m,n)s_{m,n}$ for some non-negative integers $m$ and $n$. We first notice that if $m$ and $n$ are positive then by looking at \eqref{caso m>0 n>0 verison 1} or \eqref{caso m>0 n>0 verison 2} we get $P_{x,y}(q)=1$, since otherwise $Z_{x,y}^3\neq \emptyset$. It follows that $m=0$ or $n=0$. We will prove that one of the following six cases occur:
\begin{equation} \label{path cases}
    \begin{array}{l}
     m=0, \, n=0, \mbox{ and } x=s_0;\\
     m=0, \, n=0, \mbox{ and }x=\operatorname{id};\\
     m>0,\, n=0, \mbox{ and } x=\rho(\theta(m-1,0));\\
    m>0,\, n=0, \mbox{ and } x=\rho(x_{2m});\\
      m=0,\, n>0, \mbox{ and } x=\rho^2(\theta(0,n-1));\\
  m=0,\, n>0, \mbox{ and } x=\rho^2(\sigma(x_{2n})) .
    \end{array}
\end{equation}
If $m=n=0$ then \eqref{eq theta dos 00} yields the first two cases in \eqref{path cases}. Let us now assume that $m>0$ and $n=0$. In this case we look at \eqref{caso m>0 n=0 verison 1} or \eqref{caso m>0 n=0 verison 2} in order to conclude that the elements $x\leq y$ such that $Z_{x,y}^3= \emptyset$ and  $P_{x,y}(q)\neq 1$ are those satisfying 
\begin{equation}\label{condition pat elements}
 \qquad x \not \leq  s_0\theta(m-1,0), \qquad x \not \leq \rho^2( \theta(m-1,0))s_{m,0} , \quad \mbox{ and } \qquad x \leq  \rho(\theta(m-1,0)).
\end{equation}
A direct computation shows that the only elements satisfying \eqref{condition pat elements} are $\rho(\theta(m-1,0))$ and $\rho(x_{2m})$. This gives us the third and fourth cases in \eqref{path cases}. The case when $m=0$ and $n>0$ gives us the fifth and sixth cases in \eqref{path cases} and this is treated in a similar fashion.  Finally, claims (1), (2) and (3) are clear from \eqref{path cases}, \eqref{condition pat elements}, and Proposition \ref{proposition KL theta dos}.
\end{proof}

\section{The proof}\label{proof}

We have now all the tools needed  for proving the combinatorial invariance conjecture for $\tilde{A}_2$. For the reader's convenience let us comment on our strategy. Given two isomorphic intervals $[x,y]$ and $[x',y']$ we are going to prove that their corresponding Kazhdan-Lusztig polynomials coincide by splitting the argument with respect to the location of $y$ and $y'$ in the different regions $X$, $\Theta $, $\Theta_1$, and $\Theta_2$.\\
We first treat the cases where $y$ and $y'$ are located in different regions. This is the content of  Lemma  \ref{lemdifferent regions}. These are the ``easy cases''. Indeed, we will see in the proof of Lemma \ref{lemdifferent regions} that when $y$ and $y'$ are located in different regions, a poset isomorphism between $[x,y]$ and $[x',y']$ is quite pathological and that if such an isomorphism exists, then the intervals are of small length or the corresponding Kazhdan-Lusztig polynomial is very simple ($1$ or $1+q$).

\begin{lem}\label{lemdifferent regions}
Suppose  $\phi\colon [x,y]\longrightarrow [x',y']$ is an isomorphism of posets. If $y$ and $y'$ are in different regions, then $P_{x,y}(q)=P_{x',y'}(q)$. 
\end{lem}
\begin{proof}
We split the proof in several cases in accordance with the location of $y$ and $y'$.  \\

\noindent
\textbf{Case A. } $y\in \Theta$ and $y'\in\Theta_1\cup X$. We can assume $y=\theta(m,n)$. By \eqref{KL thetas} we know that $Z_{x,y}^3=\emptyset$. By Lemma \ref{lemma Z to Z} we conclude that $Z_{x',y'}^3=\emptyset$ as well. Therefore, Lemma \ref{lemma Theta 2 with 1 in Z} implies $P_{x',y'}(q)=1$. By Proposition \ref{prop trivialKL} we get $P_{x,y}(q)=1$.\\

\noindent
\textbf{Case B. } $y\in \Theta$ and $y'\in\Theta_2$. We can assume $y=\theta(m,n)$. By \eqref{KL thetas} we know that $Z_{x,y}^3=\emptyset$. By Lemma \ref{lemma Z to Z} we conclude that $Z_{x',y'}^3=\emptyset$ as well. \\
By Proposition \ref{prop trivialKL}, we can assume that $P_{x',y'}(q)\neq 1$. Since $Z_{x',y'}^3=\emptyset$ we obtain via Lemma \ref{wrong lemma dos} that $P_{x',y'}(q) =1+q$ and $|Z_{x',y'}^4|=1$. Furthermore, $\ell (y')-\ell(x') = 4$ or $\ell (y')-\ell(x') =5$ (of course, this implies that $\ell (y)-\ell(x) = 4$ or $\ell (y)-\ell(x) =5$).\\
By Lemma \ref{lemma Z to Z} we have $|Z_{x,y}^4|=1$ and by looking at the formula in \eqref{KL thetas} we obtain $Z_{x,y}^4=\{\theta(m-1,n-1)\}$. Finally, the length constraints,  $Z_{x,y}^4=\{\theta(m-1,n-1)\}$ and \eqref{KL thetas} yield $P_{x,y}(q)=1+q$.\\

\noindent
\textbf{Case C. } $y\in \Theta_1$ and $y'\in X$.
We can assume $y=\theta(m,n)s_{m,n}$ and $y'=x_k$. Proposition \ref{proposition KL thetas uno} implies $|Z^{3}_{x,y}|\leq 2$. If $|Z^{3}_{x,y}|\leq 1$ then we can combine Lemma \ref{lemma Z to Z}, Lemma \ref{wrong lemma uno}, and Lemma \ref{lemma Theta 2 with 1 in Z} in order to obtain $P_{x,y}(q)=P_{x',y'}(q)$. Therefore, we can assume that $|Z_{x,y}^3|=2$.\\
 By  Proposition \ref{proposition KL thetas uno} we conclude that $m>0$, $n>0$ and that 
 \begin{equation*}
Z_{x,y}^3= \{ \theta(m-1,n),\theta(m,n-1)  \}.
\end{equation*} 
In particular, $x\leq \theta(m-1,n)$ and $x\leq \theta(m,n-1)$. Lemma \ref{lemma intersecccion} implies $x\leq \theta(m-1,n-1)s_{m-1,n-1} $. On the other hand,  Proposition \ref{proposition KL thetas uno} yields
\begin{equation}
    P_{\theta(m-1,n-1)s_{m-1,n-1},y}(q) = 1+ 2q.
\end{equation}
We stress that $\ell (y) -\ell ( \theta(m-1,n-1)s_{m-1,n-1} )=4$. Therefore, we can apply  Proposition \ref{proposition menor que cuatro} in order to obtain
\begin{equation}
  P_{\phi(\theta(m-1,n-1)s_{m-1,n-1}),y'}(q) =   P_{\theta(m-1,n-1)s_{m-1,n-1},y}(q)=1+2q.  
\end{equation}
This is impossible by Remark \ref{remark 1 y 1+q en X}. This finishes the proof of this case.\\

% \noindent
% \textbf{Case C. } $y\in \Theta_1$ and $y'\in X$.
% We can assume $y=\theta(m,n)s_{m,n}$ and $y'=x_k$. An inspection of \eqref{KL x} reveals that $|Z^{3}_{x',y'}|\leq 2$. If $|Z^{3}_{x',y'}|\leq 1$ then we can combine Lemma \ref{lemma Z to Z}, Lemma \ref{wrong lemma uno}, and Lemma \ref{lemma Theta 2 with 1 in Z} in order to obtain $P_{x,y}(q)=P_{x',y'}(q)$. Therefore, we can assume that $|Z_{x',y'}^3|=2$.\\
% By Lemma \ref{lemma Z to Z}, we have $|Z_{x,y}^3|=2$. By  Proposition \ref{proposition KL thetas uno} we conclude that $m>0$, $n>0$ and that $Z_{x,y}^3= \{ \theta(m-1,n),\theta(m,n-1)  \}$. In particular, $x\leq \theta(m-1,n)$ and $x\leq \theta(m,n-1)$. Lemma \ref{lemma intersecccion} implies $x\leq \theta(m-1,n-1)s_{m,n} $. On the other hand,  Proposition \ref{proposition KL thetas uno} yields
% \begin{equation}
%     P_{\theta(m-1,n-1)s_{m,n},y}(q) = 1+ 2q.
% \end{equation}
% We stress that $\ell (y) -\ell ( \theta(m-1,n-1)s_{m,n} )=4$. Therefore, we can apply  Proposition \ref{proposition menor que cuatro} in order to obtain
% \begin{equation}
%   P_{\phi(\theta(m-1,n-1)s_{m,n}),y'}(q) =   P_{\theta(m-1,n-1)s_{m,n},y}(q)=1+2q.  
% \end{equation}
% This is impossible by Remark \ref{remark 1 y 1+q en X}. This finishes the proof of this case.\\

\noindent
\textbf{Case D.} $y\in \Theta_2$ and $y'\in X$. We can assume $y=s_0\theta(m,n)s_{m,n}$ and $y'=x_k$. An inspection of \eqref{KL x} reveals that $|Z^{3}_{x',y'}|\leq 2$. If $|Z_{x',y'}^3|=1$ we can combine Lemma \ref{lemma Z to Z} and Lemma \ref{wrong lemma uno} to conclude that $P_{x,y}(q)=P_{x',y'}(q)=1+q$. Suppose now that $Z_{x',y'}^3=\emptyset$. By Lemma  \ref{lemma Theta 2 with 1 in Z} we have $P_{x',y'}(q)=1$. Then, Proposition \ref{prop trivialKL} implies $P_{x,y}(q)=1$, showing the equality of KL-polynomials in this case as well. Therefore, we can assume that $|Z_{x',y'}^3|=2$.\\
By Proposition \ref{propo x}, we must have $k$  even and greater than or equal to $6$. Furthermore,  $Z_{x',y'}^3=\{x_{k-3}, s_1s_0x_{k-5}\} $. We claim that this is impossible. On the one hand, we have \begin{equation} \label{Joins the equis son four}
  \joins{x_{k-3}}{s_1s_0x_{k-5}}{x'}{y'}{2}   =  \left\{ x_{k-1}, \rho (x_{k-1}),  \theta \left(\frac{k}{2}-2,0 \right), \rho^2  \theta \left(\frac{k}{2}-2,0 \right)  \right\}, 
\end{equation}
and therefore, $  |\joins{x_{k-3}}{s_1s_0x_{k-5}}{x'}{y'}{2}  |=4$. \\
On the other hand, by Proposition \ref{proposition KL theta dos} and Lemma \ref{lemma Z to Z} we have $Z_{x,y}^3=\{z_i,z_j\}$ for some $1\leq i\neq j \leq 4$,  where $z_i$ is defined as in Lemma \ref{Joins of theta two}. However,  Lemma \ref{Joins of theta two} implies that $| \joins{z_i}{z_j}{x}{y}{2} |\leq 3$. This contradicts Lemma \ref{lemma Join to Join} and finishes the proof of this case.\\

% By Lemma \ref{lemma Z to Z} we have $|Z_{x,y}^3|=2$.  By looking at the formulas in Proposition \ref{proposition KL theta dos} we know that 
% \begin{equation}
%     Z_{x,y}^3  = \{   z_1,z_2,z_3,z_4   \} \cap [x,y],
% \end{equation}
% where $z_i$ is defined as in Lemma \ref{Joins of theta two}. If  $Z_{x,y}^3=\{z_i,z_j\}$ then  Lemma \ref{Joins of theta two} implies that $|F_{z_i,z_j}^{[x,y]}(2)|=2$ or $|F_{z_i,z_j}^{[x,y]}(2)|=3$. However, a direct computation shows that
% \begin{equation} \label{Joins the equis son four}
%     F_{x_{k-3},s_1s_0x_{k-5}}^{[x',y']}(2) =  \left\{ x_{k-1}, \rho (x_{k-1}),  \theta \left(\frac{k}{2}-2,0 \right), \rho^2  \theta \left(\frac{k}{2}-2,0 \right)  \right\}, 
% \end{equation}
% and therefore, $  | F_{x_{k-3},s_1s_0x_{k-5}}^{[x',y']}(2)|=4$. This contradicts Lemma \ref{lemma Join to Join} and finishes the proof.
\noindent
\textbf{Case E.} $y\in \Theta_1$ and $y'\in \Theta_2$. We can assume $y=\theta(m,n)s_{m,n}$ and $y'=s_0\theta(m',n')s_{m',n'}$. By Proposition \ref{proposition KL thetas uno} we know that $|Z_{x,y}^3|\leq 2$. If $|Z_{x,y}^3|\leq 1$ we can argue as in the proof of Case D in order to conclude that $P_{x,y}(q)=P_{x',y'}(q)$. Therefore, we can assume that $|Z_{x,y}^3|=2$. \\
By Proposition \ref{proposition KL thetas uno}, we have  $Z_{x,y}^3=\{ \theta(m-1,n) , \theta(m,n-1)   \}$. Arguing as in the proof of Case C, we obtain that $\theta(m-1,n-1)s_{m-1,n-1} \in [x,y]$ and that 
\begin{equation}
    P_{\phi(\theta(m-1,n-1)s_{m-1,n-1}),y' } (q)=1+2q.
\end{equation}
 We claim that $P_{x',y'}(q) =1+q$. We remark that by \eqref{q monotonicity}  this is a contradiction that would rule out the very existence of this case, thus proving the lemma. \\
To see the above claim we consider the following elements:
\begin{align} \label{zetas primes}
z_1'&\coloneqq s_0\theta(m',n'-1); 
&z_2'&\coloneqq s_0\theta(m'-1,n'); \\
z_3'&\coloneqq \rho^2(\theta(m'-1,n'))s_{m',n'}; 
&z_4'&\coloneqq \rho(\theta(m',n'-1))s_{m',n'}.
\end{align}
These elements are the ``same'' as the ones defined in \eqref{zetas} but with the role of $m$ and $n$ replaced by $m'$ and $n'$, respectively. Proposition \ref{proposition KL theta dos} and Lemma \ref{lemma Z to Z} imply that $Z_{x',y'}^3=\{z_i',z_j'\}$ for some $1\leq i\neq j\leq 4$. On the other hand, a direct computation shows that 
  \begin{equation} \label{Joins de los dos tipos de theta s}
\joins{\theta(m-1,n)}{\theta(m,n-1)}{x}{y}{2}   = \{  \theta(m,n) , \theta(m+1,n-1), \theta(m-1,n+1)  \}   , 
 \end{equation}
and therefore $|\joins{\theta(m-1,n)}{\theta(m,n-1)}{x}{y}{2} |=3$.  Thus, by combining  Lemma \ref{lemma Join to Join}, Lemma \ref{Joins of theta two} and Lemma \ref{lemma Z to Z} we obtain that $Z_{x',y'}^3 = \{z_1',z_2' \}$ or $ Z_{x',y'}^3= \{  z_3',z_4'\}$. If $Z_{x',y'}^3= \{z_1',z_2' \}$ (resp. $Z_{x',y'}^3 = \{z_3',z_4' \}$) then using \eqref{caso m>0 n>0 verison 1} (resp. \eqref{caso m>0 n>0 verison 2}) we conclude that $P_{x',y'}(q)=1+q$. This proves our claim and finishes the proof of the lemma. 
\end{proof}

\begin{thm}\label{main theorem} In an affine Weyl group of type $\tilde{A}_2$, if $[x,y]\simeq [x',y']$ then $P_{x,y}(q)=P_{x',y'}(q)$.
\end{thm}

\begin{proof}
Throughout the proof we fix an arbitrary  poset isomorphism $\phi\colon [x,y]\longrightarrow [x',y']$.  We also assume by induction that the theorem holds for intervals of length strictly less than $\ell (y) -\ell (x)$. By Lemma \ref{lemdifferent regions} we only need to consider the case when $y$ and $y'$  belong to the same region. By Proposition \ref{proposition menor que cuatro} we can assume $\ell(y)-\ell(x)\geq 5$.  \\

\noindent
\textbf{Case A. }$y,y'\in \Theta$.  We can assume $y=\theta(m,n)$ and $y'=\theta(m',n')$.\\
Let us suppose that $m=0$ or $n=0$. Using Proposition \ref{formulas big region}   we conclude that $P_{x,y}(q)=1$ and the result follows by  Proposition \ref{prop trivialKL}. The case $m'=0$ or $n'=0$ is symmetric.\\
By the previous paragraph, we can now assume that $m,n,m',n'>0$. Let $z_0=\theta(m-1,n-1)$ and $z_0'=\theta(m'-1,n'-1)$. We can  assume that $z_0\in [x,y]$ and $z_0'\in [x',y']$. Otherwise, one of the polynomials $P_{x,y}(q)$ or $P_{x',y'}(q)$ is $1$, and Proposition \ref{prop trivialKL} proves the theorem in this case as well. We have $Z_{x,y}^4=\{z_0\}$ and $Z_{x',y'}^4=\{ z_0' \}$. Therefore, Lemma \ref{lemma Z to Z} implies $\phi (z_0) =z_0'$. Finally, we obtain 
\begin{equation}
      P_{x,y}(q) = 1 +q P_{x,z_0}(q)  = 1 +q P_{x',z_0'}(q)  = P_{x',y'}(q).
\end{equation}
The first and third equalities follow by taking the coefficient of $\HH_x$ when we expand both sides of \eqref{decomposition Theta equal N plus Theta} in terms of the standard basis of $\mathcal{H}$ and using  \eqref{passage from v to q} to pass to the $q$-version. The second equality follows by induction, since $\ell(y)-\ell(x)=\ell(z_0)-\ell(x)+4 > \ell(z_0)-\ell(x)$.  This finishes the proof in this case. \\

We will work out the details of the first and third equalities: They follow directly from \eqref{KL thetas}. Since this is the first case of the proof, we will give full details here. We will show the first one, the third one being analogous. By \eqref{KL thetas} both applied to $\theta(m,n)$ and $\theta(m-1,n-1)$ we obtain the relation
\begin{equation*}
v^2\undH_{\theta(m-1,n-1)}+\NN_{\theta(m,n)}=\undH_{\theta(m,n)}.
\end{equation*}
Let us expand all in terms of the standard basis
\begin{equation*}
v^2\sum_{x\leq z_0} h_{x,z_0}(v)\HH_x+\sum_{x\leq y} v^{\ell(y)-\ell(x)}\HH_x=\sum_{x\leq y} h_{x,y}(v)\HH_x.
\end{equation*}
Let us compare at both sides the coefficient of $\HH_x$.
\begin{equation*}
v^2h_{x,z_0}(v)+ v^{\ell(y)-\ell(x)} = h_{x,y}(v).
\end{equation*}
By \eqref{passage from v to q} we have
\begin{equation*}
v^2v^{\ell(z_0)-\ell(x)}P_{x,z_0}(v^{-2})+ v^{\ell(y)-\ell(x)} = v^{\ell(y)-\ell{(x)}} P_{x,y}(v^{-2}).
\end{equation*}
Therefore,
\begin{equation*}
P_{x,y}(v^{-2}) = 1 +v^{-2} P_{x,z_0}(v^{-2}).
\end{equation*} 
This finishes the proof in this case. \\

\noindent
\textbf{Case B. }$y,y'\in \Theta_1$. We can assume  $y =\theta(m,n)s_{m,n}$ and $y'=\theta(m',n')s_{m',n'} $.\\
By Proposition \ref{proposition KL thetas uno} we have $|Z_{x,y}^3|\leq 2$. If $|Z_{x,y}^3|\leq 1$ then Lemma \ref{lemma Z to Z}, Lemma \ref{wrong lemma uno} and Lemma \ref{lemma Theta 2 with 1 in Z}  imply $P_{x,y}(q)=P_{x',y'}(q)$. Therefore, we can assume that $|Z_{x,y}^3|=|Z_{x',y'}^3|=2$. By Proposition \ref{proposition KL thetas uno} we must have $m,n,m',n'>0 $,  $Z_{x,y}^3=\{u_1,u_2\}$ and $Z_{x',y'}^3=\{u_1',u_2'\}$, where 
\begin{equation}
  u_1= \theta(m-1,n), \quad  u_2  = \theta(m,n-1), \quad   u_1' = \theta(m'-1,n'), \quad  u_2'  =  \theta(m',n'-1).   
\end{equation}
Finally, we have  
\begin{align*}
     P_{x,y}(q) & =1 +q(P_{x,u_1}(q) + P_{x,u_2}(q)) \\
                & =1 +q(P_{x',\phi (u_1)}(q) + P_{x',\phi(u_2)}(q))\\  
                & =1 +q(P_{x',u_1'}(q) + P_{x',u_2'}(q))\\
                & =  P_{x',y'}(q),
\end{align*}  
where the first and fourth equalities follow from \eqref{Kl thetas uno}, the second one by our inductive hypothesis, and the third one is a consequence of Lemma \ref{lemma Z to Z}.\\

\noindent
\textbf{Case C. }$y,y'\in X$.  This case follows by a combination of Lemma \ref{lemma Z to Z} and Remark \ref{remark 1 y 1+q en X}.\\

\noindent
\textbf{Case D. }$y,y'\in \Theta_2$.  We can assume  $y =s_0\theta(m,n)s_{m,n}$ and $y'=s_0\theta(m',n')s_{m',n'} $. \\	
% \begin{align*}
%     z_1 &  = s_0\theta(m,n-1)  &   z_1' &  = s_0\theta(m',n'-1)\\
%     z_2 &  = s_0\theta(m-1,n)  &   z_2' &  = s_0\theta(m'-1,n') \\
%     z_3 &  = \rho (\theta(m-1,n) )s_{m,n} &   z_3' &  = \rho ( \theta(m'-1,n'))s_{m',n'} \\
%     z_4 &  = \rho^2 ( (\theta(m,n-1) )s_{m,n} &   z_4' &  = \rho^2 ( \theta(m',n'-1) )s_{m',n'}. 
% \end{align*}	
We consider the elements $z_i$ and $z_j'$ defined in \eqref{zetas} and \eqref{zetas primes}, respectively. By Proposition \ref{proposition KL theta dos} we have
\begin{equation}
    Z_{x,y}^3 = \{ z_1,z_2,z_3,z_4  \}\cap [x,y] \quad  \mbox{ and }   \quad   Z_{x',y'}^3 = \{ z_1',z_2',z_3',z_4'  \}\cap [x',y']. 
\end{equation}
We split the proof in five cases in accordance with the cardinality of $Z_{x,y}^3$ (which by Lemma \ref{lemma Z to Z} coincides with $|Z_{x',y'}^3|$).\\

\noindent
\textbf{Case D0. } $|Z_{x,y}^3|=0$. By Proposition \ref{prop trivialKL} we can assume that $P_{x,y}(q)\neq 1$ and $P_{x',y'}(q)\neq 1$. Then, Lemma \ref{wrong lemma dos} implies that $P_{x,y}(q)=P_{x',y'}(q)=1+q$. \\

\noindent
\textbf{Case D1. } $|Z_{x,y}^3|=1$. In this case we have $P_{x,y}(q)=P_{x',y'}(q)=1+q$ by Lemma \ref{wrong lemma uno}.\\

\noindent
\textbf{Case D2. } $|Z_{x,y}^3|=2$. Let us suppose $Z_{x,y}^3=\{z_1,z_2\}$. By \eqref{caso m>0 n>0 verison 1} we obtain $P_{x,y}(q)=1+q$. However, if we use \eqref{caso m>0 n>0 verison 2} to compute $P_{x,y}(q)$ then we get
 \begin{equation}
     P_{x,y}(q) =1 +2q + \mbox{ higher degree terms}.
 \end{equation}
This contradiction rules out the existence of the case $Z_{x,y}^3=\{z_1,z_2\}$. By the same reasons, we also discard the case $Z_{x,y}^3=\{z_3,z_4\}$.\\
This leaves us with four cases to be checked: $Z_{x,y}^3=\{z_1,z_3\}$, $Z_{x,y}^3=\{z_1,z_4\}$, $Z_{x,y}^3=\{z_2,z_3\}$ or $Z_{x,y}^3=\{z_2,z_4\}$.\\
Let us  assume that $Z_{x,y}^3=\{z_1,z_3\}$. We have  $m>0$ and $n>0$. By \eqref{caso m>0 n>0 verison 1}-\eqref{caso m>0 n>0 verison 2} we obtain
\begin{equation}\label{case D2 A}
    P_{x,y}(q)=1+q+qP_{x,z_3}(q).
\end{equation}
By the same argument as before, $Z_{x',y'}^3=\{z_1',z_3'\}$, $Z_{x',y'}^3=\{z_1',z_4'\}$, $Z_{x',y'}^3=\{z_2',z_3'\}$ or $Z_{x',y'}^3=\{z_2',z_4'\}$. Depending on the values of $m'$ and $n'$, we use either \eqref{caso m>0 n=0 verison 1}-\eqref{caso m>0 n=0 verison 2}, \eqref{caso m=0 n>0 verison 1}-\eqref{caso m=0 n>0 verison 2}, or \eqref{caso m>0 n>0 verison 1}-\eqref{caso m>0 n>0 verison 2} to conclude that 
\begin{equation} \label{case D2 B}
        P_{x',y'}(q)=1+q+qP_{x',z_i'}(q),
\end{equation}
where $z_i'=\phi (z_3)$. Using our inductive hypothesis we get $P_{x,z_3}(q)=P_{x',z_i'}(q)$. Therefore, \eqref{case D2 A} and \eqref{case D2 B} allow us to conclude $ P_{x,y}(q)= P_{x',y'}(q)$. The remaining three cases  are treated similarly. Indeed, the only difference that may arise in those cases is that the role of \eqref{caso m>0 n>0 verison 1}-\eqref{caso m>0 n>0 verison 2} might have to be replaced by \eqref{caso m>0 n=0 verison 1}-\eqref{caso m>0 n=0 verison 2} or \eqref{caso m=0 n>0 verison 1}-\eqref{caso m=0 n>0 verison 2}. \\

\noindent
\textbf{Case D3. } $|Z_{x,y}^3|=3$. We claim that this case is impossible. Suppose  that $Z_{x,y}^3=\{z_1,z_2,z_3\}$.  By \eqref{caso m>0 n>0 verison 1} we get
 \begin{equation}
     P_{x,y}(q) =1 +2q + \mbox{ higher degree terms}.
 \end{equation}
On the other hand, using \eqref{caso m>0 n>0 verison 2} we get 
 \begin{equation}
     P_{x,y}(q) =1 +3q + \mbox{ higher degree terms}.
 \end{equation}
This contradiction rules out the existence of this case.  The remaining three cases ($Z_{x,y}^3=\{z_1,z_2,z_4\}$, $Z_{x,y}^3=\{z_1,z_3,z_4\}$ and $Z_{x,y}^3=\{z_2,z_3,z_4\}$) are treated in a similar fashion.\\

\noindent
\textbf{Case D4. } $|Z_{x,y}^3|=4$. We have $Z_{x,y}^3=\{ z_1,z_2,z_3,z_4\}$ and $Z_{x',y'}^3=\{ z_1',z_2',z_3',z_4'\}$. By combining Lemma \ref{lemma Join to Join}, Lemma \ref{Joins of theta two} and Lemma \ref{lemma Z to Z} we obtain $\phi(\{z_1,z_2 \}) =\{z_1',z_2'\} $ or $\phi(\{z_1,z_2 \}) =\{z_3',z_4'\} $. Suppose we are in the latter case (the former case being similar). Then, we have
\begin{align*}
    P_{x,y}(q)  &  = 1+q +q(P_{x,z_1}(q)+P_{x,z_2}(q)) \\
                &  = 1+q +q(P_{x',z_3'}(q) + P_{x',z_4'}(q)) \\
                &  = P_{x',y'}(q),
\end{align*}
where the first equality follows by \eqref{caso m>0 n>0 verison 2}, the second one is a consequence of our inductive hypothesis, and the third equality follows by \eqref{caso m>0 n>0 verison 1}. This is the end of the proof.
\end{proof}

\section{Acknowledgements}

We would like to thank the anonymous referees for their careful reading of the paper, their helpful corrections, and for their interesting comments that allow us to greatly improved the exposition of the paper.

% To see this we first notice that 
% \begin{equation}
%     \theta(m-1,n)_\downarrow   \cap  \theta(m,n-1)_\downarrow =  \theta(m-1,n-1)s_\downarrow.
% \end{equation}
% Therefore, by   the monotonicity of $L$, $L\leqH R$  reduces to prove

%%%%%%%%%%%%%%%%%%%%%%%%%%%%%%%%%%%%%%%%%%%%%%%%%%%%%%%%%%%%%%%%%%%%%%%%%%%%%%%%%%%%%%%%%%%%%%%%%%%%%%%%%%%%%%%%%%%%%%%%%%%%%%%%%%%%%%%%%%%%%%%%%%%%%%%%%%%%%%%%%%%%%%%%%%%%%%%%%%%%%%%%%%%%%%%%%%%%%%%%%%%%%%%%%%%%%%%%%%%%%%%%%%%%%%%	

\bibliographystyle{plain}

\begin{thebibliography}{BCM06}

\bibitem[Ara02]{Ara02}
Alberto Arabia.
\newblock
  \href{https://webusers.imj-prg.fr/~alberto.arabia/math/intersection.pdf}{{I}ntroduction
  {\`a} l'homologie d'intersection}, 2002.
\newblock pr{\'e}publication.

\bibitem[BB05]{BB05}
Anders Bj\"{o}rner and Francesco Brenti.
\newblock {\em Combinatorics of {C}oxeter groups}, volume 231 of {\em Graduate
  Texts in Mathematics}.
\newblock Springer, New York, 2005.

\bibitem[BBP21]{batistelli2021kazhdan}
Karina Batistelli, Aram Bingham, and David Plaza.
\newblock Kazhdan-{L}usztig polynomials for $\tilde{B}_2$.
\newblock {\em arXiv preprint: 2102.01278v2}, 2021.

\bibitem[BCM06]{BCM06}
Francesco Brenti, Fabrizio Caselli, and Mario Marietti.
\newblock Special matchings and {K}azhdan-{L}usztig polynomials.
\newblock {\em Adv. Math.}, 202(2):555--601, 2006.

\bibitem[BM01]{braden2001moment}
Tom Braden and Robert MacPherson.
\newblock From moment graphs to intersection cohomology.
\newblock {\em Mathematische Annalen}, 321(3):533--551, 2001.

\bibitem[Bre94]{Br94}
Francesco Brenti.
\newblock A combinatorial formula for {K}azhdan-{L}usztig polynomials.
\newblock {\em Invent. Math.}, 118(2):371--394, 1994.

\bibitem[Bre97]{Br97}
Francesco Brenti.
\newblock Combinatorial properties of the {K}azhdan-{L}usztig {$R$}-polynomials
  for {$S_n$}.
\newblock {\em Adv. Math.}, 126(1):21--51, 1997.

\bibitem[Bre02]{Bre02}
Francesco Brenti.
\newblock Kazhdan-{L}usztig polynomials: History problems, and combinatorial
  invariance.
\newblock {\em S\'{e}m. Lothar. Combin.}, 49:Art. B49b, 30, 2002.

\bibitem[Bre04]{Br04}
Francesco Brenti.
\newblock The intersection cohomology of {S}chubert varieties is a
  combinatorial invariant.
\newblock {\em European J. Combin.}, 25(8):1151--1167, 2004.

\bibitem[Car94]{Carr91}
James~B. Carrell.
\newblock The {B}ruhat graph of a {C}oxeter group, a conjecture of {D}eodhar,
  and rational smoothness of {S}chubert varieties.
\newblock In {\em Algebraic groups and their generalizations: classical methods
  ({U}niversity {P}ark, {PA}, 1991)}, volume~56 of {\em Proc. Sympos. Pure
  Math.}, pages 53--61. Amer. Math. Soc., Providence, RI, 1994.

\bibitem[dC03]{Du03}
Fokko du~Cloux.
\newblock Rigidity of {S}chubert closures and invariance of {K}azhdan-{L}usztig
  polynomials.
\newblock {\em Adv. Math.}, 180(1):146--175, 2003.

\bibitem[Dye87]{dyer1987hecke}
Matthew Dyer.
\newblock {\em Hecke algebras and reflections in Coxeter groups}.
\newblock PhD thesis, University of Sydney Department of Mathematics, 1987.

\bibitem[Dye91]{Dy91}
Matthew Dyer.
\newblock On the ``{B}ruhat graph'' of a {C}oxeter system.
\newblock {\em Compositio Math.}, 78(2):185--191, 1991.

\bibitem[EW14]{EWhodge}
Ben Elias and Geordie Williamson.
\newblock The {H}odge theory of {S}oergel bimodules.
\newblock {\em Ann. of Math. (2)}, 180(3):1089--1136, 2014.

\bibitem[Inc06]{incitticombinatorial}
Federico Incitti.
\newblock On the combinatorial invariance of {K}azhdan--{L}usztig polynomials.
\newblock {\em J. Comb. Theory Ser. A.}, 113:1332--1350, 2006.

\bibitem[Inc07]{incitti2007more}
Federico Incitti.
\newblock More on the combinatorial invariance of {K}azhdan--{L}usztig
  polynomials.
\newblock {\em J. Comb. Theory Ser. A.}, 114(3):461--482, 2007.

\bibitem[KL79]{KL79}
David Kazhdan and George Lusztig.
\newblock Representations of {C}oxeter groups and {H}ecke algebras.
\newblock {\em Invent. Math.}, 53(2):165--184, 1979.

\bibitem[Kum02]{kum02}
Shrawan Kumar.
\newblock {\em Kac-{M}oody groups, their flag varieties and representation
  theory}, volume 204 of {\em Progress in Mathematics}.
\newblock Birkh\"{a}user Boston, Inc., Boston, MA, 2002.

\bibitem[LP20]{libedinsky2020affine}
Nicolas Libedinsky and Leonardo Patimo.
\newblock On the affine {H}ecke category for {$SL_3$}.
\newblock {\em arXiv preprint: 2005.02647v3}, 2020.

\bibitem[LPP21]{LPP21}
Nicolas Libedinsky, Leonardo Patimo, and David Plaza.
\newblock Pre-canonical bases on affine {H}ecke algebras.
\newblock {\em arXiv preprint: 2103.06903v2}, 2021.

\bibitem[Mar06]{Ma06}
Mario Marietti.
\newblock Boolean elements in {K}azhdan-{L}usztig theory.
\newblock {\em J. Algebra}, 295(1):1--26, 2006.

\bibitem[Mar18]{Ma18}
Mario Marietti.
\newblock The combinatorial invariance conjecture for parabolic
  {K}azhdan-{L}usztig polynomials of lower intervals.
\newblock {\em Adv. Math.}, 335:180--210, 2018.

\bibitem[Pat19]{Patimo2021}
Leonardo Patimo.
\newblock {A Combinatorial Formula for the Coefficient of q in Kazhdan--Lusztig
  Polynomials}.
\newblock {\em International Mathematics Research Notices}, 2021(5):3203--3223,
  10 2019.

\bibitem[Pla17]{plaza2017graded}
David Plaza.
\newblock Graded cellularity and the monotonicity conjecture.
\newblock {\em J. Algebra}, 473:324--351, 2017.

\bibitem[Soe97]{soergel1997kazhdan}
Wolfgang Soergel.
\newblock Kazhdan-{L}usztig polynomials and a combinatoric for tilting modules.
\newblock {\em Representation Theory of the AMS}, 1(6):83--114, 1997.

\end{thebibliography}

\vspace{1cm}
\sc
\noindent
gbur0996@uni.sydney.edu.au, The University of Sydney, Australia. \newline
nlibedinsky@u.uchile.cl, Universidad de Chile, Chile. \newline
dplaza@inst-mat.utalca.cl, Universidad de Talca, Chile. 

\end{document}